\documentclass[12pt, reqno]{amsart}
\usepackage{mathrsfs}
\usepackage{amsfonts}
\usepackage[centertags]{amsmath}
\usepackage{amssymb,array}
\usepackage{amsthm}
\usepackage{stmaryrd}
\usepackage{graphicx}
\usepackage{caption}
\usepackage{ytableau}
\usepackage{enumerate,enumitem}
\usepackage{MnSymbol}
\SetLabelAlign{center}{\hfill #1\hfill}
\usepackage[textwidth=16cm, hmarginratio=1:1]{geometry}
\usepackage{tikz-cd, tikz}
\usepackage{rotating}
\usepackage[all,cmtip]{xy}
\usepackage[titletoc, title]{appendix}
\usepackage{amscd}
\usepackage[colorlinks]{hyperref}
\usepackage{float}
\usepackage{dsfont} 
\usepackage{comment}
\hypersetup{
bookmarksnumbered,
pdfstartview={FitH},
breaklinks=true,
linkcolor=blue,
urlcolor=blue,
citecolor=blue,
bookmarksdepth=2
}

\theoremstyle{plain}
   
   \newtheorem{theorem}{Theorem}[section]
   \newtheorem{proposition}[theorem]{Proposition}
   \newtheorem{prop}[theorem]{Proposition}
   \newtheorem{lemma}[theorem]{Lemma}
   \newtheorem{corollary}[theorem]{Corollary}
   \newtheorem{conjecture}[theorem]{Conjecture}
   
\theoremstyle{definition}
   \newtheorem{definition}[theorem]{Definition}
   
   \newtheorem{example}[theorem]{Example}

   \newtheorem{remark}[theorem]{Remark}

\numberwithin{equation}{section}

\newcommand{\CC}{{\mathbb {C}}}

\newcommand{\RR}{{\mathbb {R}}}
\newcommand{\ZZ}{{\mathbb {Z}}}

\newcommand{\FF}{{\mathcal{F}}}

\newcommand{\II}{{\mathcal {I}}}

\newcommand{\ch}{{\operatorname{ch}}}

\newcommand{\SSYT}{{\rm SSYT}}

\newcommand{\m}{{\bf m}}
\newcommand{\n}{{\bf n}}

\newcommand{\om}{\overline{{\bf m}}}
\newcommand{\Z}{{\rm Z}}
\newcommand{\oZ}{\overline{{\rm Z}}}

\DeclareMathOperator{\rep}{Rep}

\DeclareMathOperator{\id}{id}

\DeclareMathOperator{\sgn}{sgn}

\DeclareMathOperator{\Gr}{Gr}
\DeclareMathOperator{\pr}{\rm pr}
\DeclareMathOperator{\Irr}{Irr}

\DeclareMathOperator{\Mult}{Mult}
\DeclareMathOperator{\soc}{\rm soc}
\DeclareMathOperator{\coso}{\rm cos}

\newcommand\scalemath[2]{\scalebox{#1}{\mbox{\ensuremath{\displaystyle #2}}}}

\newlength{\mysizetiny}
\newlength{\mysizesmall}
\newlength{\mysize}
\newlength{\mysizelarge}
\begin{document}

\title[On the simplicity of tensor products]{On the simplicity of the tensor product of two simple modules of quantum affine algebras}
\author{L\'{e}a Bittmann}
\address{L\'{e}a Bittmann, Université de Strasbourg, 7 rue René Descartes, 67000 Strasbourg, France}

\author{Jian-Rong Li}
\address{Jian-Rong Li, Faculty of Mathematics, University of Vienna, Oskar-Morgenstern-Platz 1, 1090 Vienna, Austria.} 

\date{} 

\begin{abstract}
By Chari and Pressley's quantum affine Schur Weyl duality, the representation theory of type $A$ quantum affine algebras is equivalent to that of affine Hecke algebras, Which is in turn equivalent, through type theory, to the representation theory of general linear groups over a non-archimedean field. One caveat being that the equivalence only stands for representations of the quantum affine algebra $U_q(\widehat{\mathfrak{sl}_k})$, with $k$ large enough. This paper investigates the case of smaller $k$. In particular, we generalize a criterion of Lapid and Minguez for the irreducibility of the parabolic induction of two ladder representations and obtain a necessary and sufficient condition for the irreducibility of the tensor product of two snake modules, for the quantum affine algebra $U_q(\widehat{\mathfrak{sl}_k})$. The criterion also work when one module is a fundamental module at an extremity node and the other is any simple module. The question of the irreducibility of the tensor product of simple modules is important in the study of representations of quantum affine algebras, in particular with their relations to cluster algebras. 
Additionally, we translates our criterion to the setting of Grassmannian cluster algebras $\mathbb{C}[\mathrm{Gr}(k,n, \sim)]$, where it generalizes Leclerc and Zelevinsky's result that two Pl\"{u}cker coordinates are compatible if and only if they are weakly separated.
\vspace{-1cm}

\end{abstract}

\maketitle 

\setcounter{tocdepth}{1}
\tableofcontents
\vspace{-1cm}

\section{Introduction}

Let $\mathfrak{g}$ be a simple Lie algebra, $q\in\mathbb{C}^\times$ be not a root of unity and consider $U_q(\widehat{\mathfrak{g}})$, the corresponding quantum affine algebra. We are interested in the tensor category $\mathscr{C}$ of (type 1) finite-dimensional representations (these representations can also be regarded as representations of a quantum loop algebra $U_q(\mathcal{L}\mathfrak{g})$).
 These simple finite-dimensional representations were classified by Chari-Pressley \cite{CP94}, they are highest ($l$-)weight representations, with weights in a certain monoid. However, some important questions about the category $\mathscr{C}$ remain open. For example, it is not known, in general, under which conditions the tensor product $L(M)\otimes L(N)$ of two irreducible finite-dimensional representations is also irreducible. Additionally, for questions related to the cluster algebra structure of the Grothendieck ring $K_0(\mathscr{C})$ (see \cite{HL10,HL16}), it is important to know which irreducible representations $L$ are \textit{real}, i.e. its tensor square $L(M)\otimes L(M)$ is irreducible.
We focus here on the case of a type $A$ quantum affine algebra, that is $\mathfrak{g}=\mathfrak{sl}_k$.

Chari and Pressley \cite{CP96b} defined a quantum affine Schur Weyl duality functor between the category of finite-dimensional modules for the affine Hecke algebra $\hat{H}_N(q^2)$ and the category $\mathscr{C}^N_{\mathfrak{sl}_k}$ of finite-dimensional $U_q(\widehat{\mathfrak{sl}}_k)$-modules of level $N$. If $N< k$, this functor is an equivalence of categories.  

On the other hand, by \emph{type theory} (see  \cite{Bor76} \cite{Cas80}, \cite{BK98} and \cite{SS12}), or equivalentely the work of Heiermann \cite{H11}, there exists an equivalence of categories between smooth, finite-length representations of $GL_N(F)$, where $F$ is a non-archimedian field, in a simple Bernstein block, and the category of finite-dimensional right $\hat{H}_N(q^2)$-modules. 

Thus, for $k>N$, there is an equivalence of categories between smooth, finite-length representations of $GL_N(F)$ in a simple Bernstein block and the category $\mathscr{C}^N_{\mathfrak{sl}_k}$ of finite-dimensional $U_q(\widehat{\mathfrak{sl}}_k)$-modules of level $N$, and this equivalence is \textit{monoidal}: the tensor product is translated into the parabolic induction.  Therefore results of representations of general linear groups over $F$ can be applied to study representations of $U_q(\widehat{\mathfrak{sl}_k})$, when $k$ is large enough. 

For smaller values of $k$, the quantum affine Schur-Weyl duality functor is not an equivalence and this is the case we mostly focus on in the present paper. 
For all values of $k$, we construct a candidate category $\mathcal{C}_k$, via quotient and localization of the simple block category, which should be equivalent to the category $\mathscr{C}_{\mathfrak{sl}_k}$ of finite-dimensional representations of $U_q(\widehat{\mathfrak{sl}}_k)$, of all levels. We prove an isomorphism of Grothendieck rings between the two categories $\mathscr{C}_{\mathfrak{sl}_k}$ and $\mathcal{C}_k$ (see Proposition~\ref{prop: iso Rk's}). This allows us to obtain results on representations of $U_q(\widehat{\mathfrak{sl}_k})$, for any $k \ge 1$, by using modified results on representations of general linear groups over $F$.

In \cite{LM16}, Lapid-M\'{i}nguez considered a certain class of representations called \emph{ladders}, whose equivalence under the quantum affine Schur Weyl duality are studied under the name \emph{snake} modules (see \cite{MY12}). They established a combinatorial criterion for the irreducibility of the parabolic induction $\sigma \times \pi$, where $\sigma$ is a ladder representation and $\pi$ is an arbitrary irreducible representation of $GL_N(F)$. Thus, through quantum affine Schur-Weyl duality, when $k$ is large enough, Lapid and M\'{i}nguez's criterion gives a criterion of the simplicity of the tensor product of a snake module and any simple module of the quantum affine algebra $U_q(\widehat{\mathfrak{sl}_k})$.



In this paper, we add conditions to Lapid and M\'{i}nguez's criterion to obtain a criterion which works for any $k \ge 1$.

Our main result is the following.

\begin{theorem}[{Theorem \ref{thm:condition of irreducibility of tensor product}}, Proposition~\ref{prop:condition of irred for Yk-1}] \label{thm:main theorem in introduction}
Let $L(M), L(M')$ be simple $U_q(\widehat{\mathfrak{sl}_k})$-modules. Assume either $L(M')$ is a fundamental representation at an extremal node, or both $L(M)$ and $L(M')$ are snake modules. Then the tensor product $L(M)\otimes L(M')$ is irreducible if and only if the highest weights $M,M'$ satisfy two combinatorial conditions ${\rm LC}_k(\m_M,\m_{M'})$ and ${\rm LC}_k(\m_{M'}, \m_M)$.
\end{theorem}

In order to prove this result, we use our result on the isomorphism of Grothendieck rings, as well as the results of Maxim Gurevich \cite{Gur21}, who gave an algorithm to compute the decomposition of a tensor product of any two ladder representations into irreducible representations.



Additionally, Hernandez and Leclerc \cite {HL10} showed that in the case of a quantum affine algebra of type $A$, there is an isomorphism between the cluster algebra structure of the Grothendieck ring of certain subcategories $\mathscr{C}_\ell$ of $\mathscr{C}$ and some Grassmannian cluster algebras. More precisely, there is an isomorphism between $K_0(\mathscr{C}_{\ell})$ and $\mathbb{C}[\Gr(k,n,\sim)]$, where $n = k+ \ell +1$ and $\mathbb{C}[\Gr(k,n,\sim)]$ is the quotient of $\mathbb{C}[\Gr(k,n)]$ by the ideal generated by $P_{i,i+1, \ldots, i+k-1}-1$, $i \in [n-k+1]$, where $P_{i,i+1, \ldots, i+k-1}$ is a Pl\"{u}cker coordinate.

It is shown in \cite{CDFL} that the dual canonical basis of $\mathbb{C}[\Gr(k,n,\sim)]$ is parametrized by some rectangular semistandard Young tableaux. In particular, cluster monomials in $\mathbb{C}[\Gr(k,n,\sim)]$ are of the form $\ch(T)$, where $T$ is such a tableau. Two cluster variables $\ch(T)$, $\ch(T')$ are called \emph{compatible} if they appears in the same cluster, and this implies that $\ch(T)\ch(T') = \ch(T \cup T')$. It is conjectured in \cite{CDFL} that $\ch(T)\ch(T') = \ch(T \cup T')$ is equivalent to the condition that $\ch(T)$, $\ch(T')$ are compatible. 

The result of Theorem~\ref{thm:main theorem in introduction} translates to the context of the Grassmannian cluster algebra. We obtain a combinatorial criterion on the ladders tableaux $T,T'$ equivalent to the condition $\ch(T)\ch(T') = \ch(T \cup T')$ (Corollary~\ref{cor:main result tableaux}). In particular, it generalizes a result of Leclerc-Zelevinsky \cite{LZ} stating that any two Pl\"{u}cker coordinates are compatible if and only if they are weakly separated.


Note that by adopting the approach of \cite{LM20}, our results might be generalized to quantum affine superalgebras. Then $L(M)$ should be a skew representation while $L(M')$ needs to be replaced by any simple module that is a subquotient of tensor products of evaluation vector representations.

The paper is organized as follows. We start by recalling results on representations of quantum affine algebras and $p$-adic groups in Section~\ref{sec:quantum affine algebras and p-adic groups}. Then, Section~\ref{sec:quantum affine Schur Weyl duality} is dedicated to the quantum affine Schur-Weyl duality, and how these representations are related. We establish the isomorphism of Grothendieck rings that will be used later on. In Section~\ref{sec: QA and k-Z classification}, we establish an analog of the Zelevinsky classification \cite{Zel}, in terms of representations for the type $A$ quantum affine algebra, and for our quotient category $\mathcal{C}_k$. In Section~\ref{sec:a criterion for simplicity of tensor products}, we describe the conjectured criterion for the simplicity of the tensor product of two $U_q(\widehat{\mathfrak{sl}_k})$-modules, and obtain the main results. Section~\ref{sect_proof} contains the proofs of the two main intermediate results. In Section \ref{sec:application to Grassmannian cluster algebras}, we apply our criterion of simplicity of tensor product of two snake modules to obtain compatibility criterion of two cluster variables in a Grassmannian cluster algebra which correspond to snake modules.

\subsection*{Acknowledgements}

We would like to thank Alberto Mínguez for very helpful discussions, Maxim Gurevich for explaining to us his results about decomposition of the product of two ladder representations and Greg Warrington for his C program of computing Kazhdan-Lusztig polynomials, which is used in the computation of $q$-characters of simple modules of $U_q(\widehat{\mathfrak{sl}_k})$. The authors would like thank the reviewers for their helpful comments and suggestions.

LB was partially supported by Austrian Science Fund (FWF) Project P-31705, by the European Research Council (ERC) under the European Union's Horizon 2020 research and innovation programme under grant agreement No 948885 and by the Royal Society University Research Fellowship, JL was supported by the Austrian Science Fund (FWF): P-34602, Grant DOI: 10.55776/P34602, and PAT 9039323, Grant-DOI 10.55776/PAT9039323.

\subsection*{Notations}
For convenience of the reader, we collect key notation here. 

\begin{itemize}
\item $U_q(\widehat{\mathfrak{g}})$ the quantum affine algebra for a complex simple Lie algebra $\mathfrak{g}$, $I$ the set of vertices of its Dynkin diagram, Section \ref{subsec:quantum affine algebras}.

\item $\hat{I} = \{ (i,s) \in I\times \ZZ \mid i +s -1 \in 2\ZZ\}$, $\mathcal{P}$ the free abelian group in formal variables $Y_{i,s}^{\pm 1}$, $(i,s) \in \hat{I}$, $\mathcal{P}^+$ the submonoid of $\mathcal{P}$ generated by $Y_{i,s}$,  $(i,s) \in \hat{I}$, Section \ref{subsec:modules of quantum affine algebras}.

\item $\mathscr{C}$ the category of finite dimensional representations of $U_q(\widehat{\mathfrak{g}})$, $\mathscr{C}^\ZZ$ a subcategory of $\mathscr{C}$ and $\mathscr{R}_\mathfrak{g}$ its Grothendieck ring, $\mathscr{R}_k = \mathscr{R}_{\mathfrak{sl}_k}$, Section \ref{subsec:modules of quantum affine algebras},

\item $L(M)$ the simple $U_q(\widehat{\mathfrak{g}})$-module with highest $l$-weight $M$ and $\chi_q(M) = \chi_q(L(M))$ its $q$-character, Section \ref{subsec:modules of quantum affine algebras}.

\item $\mathcal{C}$ the category of complex, smooth representations of $GL_N(F)$ ($N=0,1,2,\ldots$) of finite length, $\mathcal{C}^{\ZZ}$ a certain subcategory of $\mathcal{C}$ and $\mathcal{C}_k$ a certain quotient of $\mathcal{C}^{\ZZ}$, $\mathcal{R}$ and $\mathcal{R}_k$ the respective Grothendieck rings of $\mathcal{C}^\ZZ$ and $\mathcal{C}_k$, Sections \ref{sect_p_adic}, \ref{sec:quantum affine Schur Weyl duality}.

\item For $\Delta = [c,d]_{\rho}$, $b(\Delta) = \rho \nu_{\rho}^c$, $e(\Delta) = \rho \nu_{\rho}^d$, $\overleftarrow{\Delta} = [c-1, d-1]_{\rho}$, ${}^-\! \Delta = [c+1,d]_{\rho}$, $\nu_\rho$ is the character $\nu_\rho(g) = |\det(g)|^{s_{\rho}}$, $s_\rho$ is a certain element in $\RR_{>0}$, Section \ref{sect_p_adic}.

\item Left aligned order $\le_b$ on segments, right aligned order $\le_e$ on segments, precede order $\prec$ on segments, $k$-precede order on segments $\prec_k$, Sections \ref{sect_p_adic}, \ref{subsec:Zelevinsky classification in Ck}.

\item ${\rm Irr}$ is the set of irreducible representations in $\mathcal{C}$ up to equivalence. ${\rm Irr}_k$ is the set of irreducible representations in $\mathcal{C}_k$ up to equivalence, Sections \ref{sect_p_adic}, \ref{sec:quantum affine Schur Weyl duality}.

\item $\Gr(k,n)$ the Grassmannian of $k$-planes in $\mathbb{C}^n$ and $\CC[\Gr(k,n)]$ its homogeneous coordinate ring; $\CC[\Gr(k,n,\sim)]$ the quotient of $\CC[\Gr(k,n)]$ by the Pl\"{u}cker coordinates with column set a consecutive interval; $P_{i_1, \ldots, i_n} \in \CC[{\rm Gr}(k, n)]$ a Pl\"ucker coordinate, Section \ref{subsec:Grassmannian cluster algebras and semistandard Young tableaux}. 

\item ${\rm SSYT}(k, [n])$ the monoid of rectangular semistandard Young tableaux with $k$ rows and with entries in $[n]$; ${\rm SSYT}(k, [n],\sim)$ the monoid of $\sim$-equivalence classes, Section \ref{subsec:Grassmannian cluster algebras and semistandard Young tableaux}. 

\item $\ch(T)$ is the dual canonical basis element in $\CC[\Gr(k,n \sim)]$ corresponding to a tableau $T \in {\rm SSYT}(k, [n],\sim)$, Section \ref{subsec:Grassmannian cluster algebras and semistandard Young tableaux}. 

\item $\rightsquigarrow$-matching, $\rightsquigarrow$-matching function, best $\rightsquigarrow$-matching, Section \ref{subsec:matching functions}.

\item $X_{{\bf m}, {\bf n}}$, $X_{{\bf m}, {\bf n}}^{(k)}$, $X_{T,T'}$, are certain sets defined in Sections \ref{subsec:LCmn, LCkmn} and \ref{subsec:condition LCTTprime}.

\item ${\rm LC}({\bf m}, {\bf n})$, ${\rm LC}_k({\bf m}, {\bf n})$, ${\rm LC}(T, T')$ are certain conditions defined in Sections \ref{subsec:LCmn, LCkmn} and \ref{subsec:condition LCTTprime}.

\item ${\rm LI}({\bf m}, {\bf n})$, ${\rm RI}({\bf m}, {\bf n})$, ${\rm LI}_k({\bf m}, {\bf n})$, ${\rm RI}_k({\bf m}, {\bf n})$ are certain conditions defined in Section \ref{subsec:LI,RI}.



\end{itemize}

\section{Quantum affine algebras and \texorpdfstring{$p$}{p}-adic groups}  \label{sec:quantum affine algebras and p-adic groups}

\subsection{Quantum affine algebras}\label{definition of quantum affine algebras}\label{subsec:quantum affine algebras}

Let $\mathfrak{g}$ be a simple Lie algebra of simply laced type, and $I=\{ 1, \ldots, \ell\}$ be the vertices of its Dynkin diagram. Let $\hat{\mathfrak{g}}$ be the associated affine Lie algebra and $C=(C_{ij})_{0\leq i,j\leq \ell}$ the Cartan matrix of $\hat{\mathfrak{g}}$. Fix $q\in\mathbb{C}^\times$ which is not a root of unity.

The $q$-numbers, $q$-factorials and $q$-binomial coefficients are defined as follows:
\begin{equation*}
[n]_q = \frac{q^n-q^{-n}}{q-q^{-1}}, \quad [n]_q!=[h]_q \cdots [1]_q, \quad {\binom{m}{n}}_q = \frac{[m]_q!}{[n]_q! [m-n]_q!}.
\end{equation*}


The \emph{quantum affine algebra} $U_q(\widehat{\mathfrak{g}})$ in Drinfeld's realization \cite{Dri2} is the $\mathbb{C}$-algebra with generators $x_i^{\pm}$, $k_i^{\pm 1}$ ($i=0, \ldots, \ell$) and relations:
\begin{equation*}
\begin{array}{cc}
k_i k_i^- = k_i^- k_i = 1, \quad k_ik_j = k_jk_i, & k_i x_j^{\pm} k_i^{-1}  = q^{\pm r_iC_{ij}} x_j^{\pm}, \\
{[}x_i^+, x_j^-] = \delta_{ij} \frac{k_i-k_i^{-1}}{q-q^{-1}}, & \sum_{r=0}^{1-C_{ij}} (-1)^r {\binom{1-C_{ij}}{r}}_{q}  (x_i^{\pm})^r x_j^{\pm} (x_i^{\pm})^{1-C_{ij}-r} = 0.
\end{array}
\end{equation*}

The algebra $U_q(\widehat{\mathfrak{g}})$ has a structure of a Hopf algebra with the comultiplication $\Delta$ and antipode $S$ given on the generators by the formulas:
\[
\begin{array}{cc}
\Delta(k_i) = k_i \otimes k_i, & S(k_i^{\pm 1}) = k_i^{\mp 1},\\
\Delta(x_i^+) = x_i^+ \otimes 1 + k_i \otimes x_i^+, & S(x_i^+) = - x_i^+ k_i,\\
\Delta(x_i^-) = x_i^- \otimes k_i^{-1} + 1 \otimes x_i^-, & S(x_i^-) = -k_i^{-1} x_i^-.
\end{array}
\]

\subsection{Finite dimensional modules of quantum affine algebras} \label{subsec:modules of quantum affine algebras}
In this section, we recall the standard facts about finite-dimensional $U_q(\widehat{\mathfrak{g}})$-modules and their $q$-characters, see \cite{CP94, CP95a, FR98}.

Let $\mathscr{C}$ be the category of finite-dimensional $U_q(\widehat{\mathfrak{g}})$-modules (of type 1). Chari-Pressley \cite{CP94} have shown that irreducible representations in this category are highest weight modules, and are  labeled by \emph{Drinfeld polynomials}. Here, we write these highest weights as \emph{dominant monomials} in $\mathbb{Z}[Y_{i,a} \mid i\in I, a\in \mathbb{C}^\times]$, $M= Y_{i_1,a_1} Y_{i_2,a_2} \cdots Y_{i_r,a_r}$. For $M$ such a dominant monomial, let $L(M)$ be the simple module of highest weight $M$. For $i\in I, a\in\mathbb{C}^\times$, the simple module $L(Y_{i,a})$ is called a \emph{fundamental} representation.

Let $\hat{I} = \{ (i,s) \in I\times \mathbb{Z} \mid i+s -1 \in 2\ZZ\}$\footnote{It corresponds to the choice of \emph{height function} $\xi_i=i-1$, see \cite{HL10}}. 
For all quantum parameters $a\in \mathbb{C}^\times$, let $\mathcal{P}^+_a$ be the dominants monomials of the form $Y_{i_1,aq^{s_1}}Y_{i_2,aq^{s_2}}\cdots Y_{i_r,aq^{s_r}}$, for $(i_j,s_j) \in \hat{I}$. For $a_1,\ldots, a_r$ such that $a_i/a_j\notin q^{2\mathbb{Z}}$ for $i\neq j$, and for $M_i\in \mathcal{P}^+_{a_i}$, the tensor product 
\[L(M_1) \otimes L(M_2) \otimes \cdots L(M_r) \]
is irreducible. Thus we can restrict ourselves to the study of the subcategory of representations in $\mathscr{C}$ whose composition factors have highest weight in one $\mathcal{P}^+_a$ (for any choice of $a$). 

From now on, for $(i,s)\in \hat{I}$, let us write $Y_{i,s}$ for $Y_{i,q^s}$. Denote by $\mathcal{P}$ the free abelian group generated by $Y_{i,s}^{\pm 1}$, $(i,s)\in \hat{I}$, and let  $\mathcal{P}^+$ the submonoid of $\mathcal{P}$ generated by $Y_{i,s}$, $(i,s)\in \hat{I}$. Let $\mathscr{R}_\mathfrak{g}$ be the Grothendieck ring of the category $\mathscr{C}^\mathbb{Z}$, it is freely generated, as a commutative ring, by the images of the fundamental modules  $L(Y_{i,s})$, for $(i,s)\in \hat{I}$. For $\mathfrak{g}=\mathfrak{sl}_k$, we use the notation $\mathscr{R}_k:= \mathscr{R}_{\mathfrak{sl}_k}$.

Let $\mathbb{Z}\mathcal{P} = \mathbb{Z}[Y_{i, s}^{\pm 1} \mid (i,s) \in \hat{I}]$ be the group ring of $\mathcal{P}$. Frenkel-Reshetikhin \cite{FR98} have introduced the \emph{$q$-character}, an injective ring morphism 
\[\chi_q : \mathscr{R}_\mathfrak{g} \to \mathbb{Z}\mathcal{P}.\]
For a $U_q(\widehat{\mathfrak{g}})$-module $V$, $\chi_q(V)$ encodes the decomposition of $V$ into common generalized eigenspaces for the action of a large commutative subalgebra of $U_q(\widehat{\mathfrak{g}})$ (the \emph{loop}-Cartan subalgebra).

\begin{example}
For $\mathfrak{g} = \mathfrak{sl}_2$, the $q$-character of the fundamental representation $L(Y_{1,0})$ is 
\[\chi_q(L(Y_{1,0})) = Y_{1,0} + Y_{1,2}^{-1}.\] 
\end{example}

\subsection{Representations of \texorpdfstring{$p$}{p}-adic \texorpdfstring{$GL_N$}{GLN}}\label{sect_p_adic}
We recall certain results about certain representations of the general linear group over a non-archimedean local field, see \cite{BZ,Zel,LM18}. 

Let $F$ be a non-archimedean local field of characteristic $p$, with a normalized absolute value $|\cdot|_F$. 
For any $N\in \mathbb{Z}_{\geq 0}$, consider the group of $F$-points $G_N=GL_N(F)$ (with the convention that $G_0$ is the trivial group). The topology on $F$ naturally equips $G$ with a local compact totally disconnected topology.

Let $\mathcal{C}(G_N)$ be the category of complex, smooth representations of $G_N$ of finite length. Denote by $\text{Irr}G_N$ the set of irreducible objects of $\mathcal{C}(G_N)$, up to equivalence. 
By abuse of notation, for any $N$, let $\nu$ be the character $\nu = |\det|_F$. For $\pi \in \mathcal{C}(G_N)$, the \emph{degree} of the representation is $\deg(\pi) = N$.

For $\pi_i \in \mathcal{C}(G_{N_i})$, $i=1,2$, denote by $\pi_1 \times \pi_2 \in \mathcal{C}(G_{N_1+N_2})$ the representation which is parabolically induced from $\pi_1 \otimes \pi_2$. The parabolic induction endows the category $\bigoplus_{N \ge 0}\mathcal{C}(G_{N})$ with the structure of a tensor category.

Denote by $\mathcal{R}^{G_N}$ (resp. $\mathcal{R}^G$) the Grothendieck ring of $\mathcal{C}(G_N)$ (resp. $\mathcal{C}=\oplus_{N \ge 0} \mathcal{C}(G_N)$). Then $\mathcal{R}^G = \oplus_{N \ge 0} \mathcal{R}^{G_N}$ is a commutative graded ring under $\times$.  Denote $\text{Irr} = \cup_{\ge 0} {\rm Irr}G_N$. 

Recall that an irreducible smooth representation $(\pi, V)$ of $G_N$ is called supercuspidal if for all proper parabolic subgroups $P \subsetneq G_N$ with Levi subgroup $M$ and all irreducible smooth representations $(\sigma, W)$ of $M$,  the representation $(\pi, V)$ is not a subrepresentation of the parabolic induction $( {\rm Ind}^{G_N}_{P}\sigma, {\rm Ind}^{G_N}_{P} W )$, see \cite[Definition 2.4.7]{Fin}. Denote by ${\rm Irr}_c \subset {\rm Irr}$ the subset of supercuspidal representations of $G_N$, $N > 0$. 

For any $\pi \in {\rm Irr}$, there exist supercuspidal representations $\rho_1, \ldots, \rho_r$, uniquely determined up to permutation, such that $\pi$ is a subrepresentation of $\rho_1 \times \cdots \times \rho_r$. Then ${\rm supp}(\pi) := \{\rho_1, \ldots, \rho_r\}$ (not multiset) is called the supercuspidal support of $\pi$ \cite{LM16}. 

For any $\rho \in {\rm Irr}_c $, there exists a unique $s_{\rho} \in \RR_{>0}$ such that, for any $\rho'\in {\rm Irr}_c $, the induced representation $\rho \times \rho'$ is reducible if and only if $\deg(\rho)=\deg(\rho')$ and either $\rho'=\rho \nu^{s_\rho}$ or $\rho = \rho' \nu^{s_\rho}$. Let $\nu_{\rho} =  \nu^{s_\rho}$,  $\overrightarrow{\rho} = \rho \nu_\rho$, $\overleftarrow{\rho}=\rho \nu_\rho^{-1}$.

A \emph{segment} is a finite nonempty subset of ${\rm Irr}_c$ of the form $\Delta = \{\rho_1, \ldots, \rho_r\}$, where $\rho_{i+1} = \overrightarrow{\rho_i}$, $i \in [r-1]$. Let $b(\Delta)=\rho_1$, $e(\Delta)=\rho_r$, $\deg(\Delta) = \sum_{i=1}^r \deg(\rho_i) = r\cdot \deg(\rho_1)$ and the \emph{length} of $\Delta$ is $\ell(\Delta) = r$. Additionally, $\overrightarrow{\Delta} = \{ \overrightarrow{\rho_1}, \ldots , \overrightarrow{\rho_r}\}$ and $\overleftarrow{\Delta} = \{ \overleftarrow{\rho_1}, \ldots , \overleftarrow{\rho_r}\}$. 

For two segments $\Delta_1,\Delta_2$, we say that $\Delta_1$ \emph{precedes} $\Delta_2$, and write $\Delta_1 \prec \Delta_2$ if and only if $b(\Delta_1) \notin \Delta_2$, $b(\overleftarrow{\Delta_2}) \in \Delta_1$ and $e(\Delta_2) \notin \Delta_1$. The segments $\Delta_1$ and $\Delta_2$ are \emph{linked} if either $\Delta_1 \prec \Delta_2$ or $\Delta_2 \prec \Delta_1$, or equivalently if the union $\Delta_1 \cup \Delta_2$ is a segment which strictly contains both segments.

For a segment $\Delta =\{\rho_1, \ldots, \rho_r\}$, let ${\rm Z}(\Delta) := {\rm soc}(\rho_1 \times \cdots \times \rho_r) \in {\rm Irr} G_{\deg \Delta}$, where ${\rm soc}(\pi)$ denotes the socle of $\pi$, i.e., the largest semisimple subrepresentation of $\pi$. We use the convention that ${\rm Z}(\emptyset)=1$, the trivial representation of $G_0$. 

A {\sl multisegment} is a formal finite sum ${\bf m} = \sum_{i=1}^m \Delta_i$ of segments. Let $\deg {\bf m} = \sum_{i=1}^m \deg(\Delta_i)$. 

 We say that the multisegment ${\bf m} = \sum_{i=1}^m \Delta_i$ is \emph{ordered} if $\Delta_i \not\prec \Delta_j$ for all $i<j$. 
For an ordered multisegment ${\bf m} = \sum_{i=1}^m \Delta_i$, let 
\begin{align*}
\zeta({\bf m}) & := {\rm Z}(\Delta_1) \times \cdots \times {\rm Z}(\Delta_m) \in \mathcal{C}(G_{\deg {\bf m}}),\\
\text{and} \quad {\rm Z}({\bf m}) & := {\rm soc}(\zeta({\bf m})) \in {\rm Irr}G_{\deg {\bf m}}
\end{align*}
 The map $\m  \mapsto \Z(\m)$ sending a multisegment to an irreducible representation
is a bijection (\emph{the Zelevinsky Classification}), see \cite{BZ,Zel}.
As a consequence of the Zelevinsky Classification, the representation $\zeta({\bf m})$ is irreducible if and only if the multisegment $\m$ is \emph{pairewise unliked}, i.e. no two segments in $\m$ are linked. 

Thus, as above, we can restrict ourselves to the study of a skeleton subcategory of $\bigoplus_{N \ge 0}\mathcal{C}(G_{N})$, without loss of information. A \emph{cuspidal line} is the equivalent class for the equivalence relation generated by $\rho \sim \overleftarrow{\rho}$. For $\rho \in \Irr_c$, the cuspidal line containing $\rho$ is $\ZZ_\rho = \{\rho \nu_{\rho}^{a}: a \in \ZZ\}$.
By the Zelevinsky Classification, if the representations $\pi_1, \pi_2, \ldots, \pi_r$ have supercuspidal support contained in different cuspidal lines, then the representation $\pi_1 \times \pi_2 \times \cdots \times \pi_r$ is irreducible. Hence in practice, it is enough to consider representations with support inside a single supercuspidal line. From now on, we fix $\rho \in \Irr_c$.

We denote by $\mathcal{C}^{\mathbb{Z}}$ the Serre subcategory of $\bigoplus_{N\geq 0}\mathcal{C}(G_N)$ of representations whose supercuspidal support is in a fixed line $\mathbb{Z}_\rho$ (from \cite[Theorem 3.8]{LM18} these categories are all equivalent so the choice of the supercuspidal representation $\rho$ does not matter), its irreducible representations are given by $\Irr^\mathbb{Z}$.


For $a\leq b \in \ZZ$, we write $[a,b]$ for the segment $\Delta = \{ \rho \nu_\rho^{a}, \rho \nu_\rho^{a+1}, \ldots, \rho \nu_\rho^{b} \}$, let $\mathrm{Seg}$ be the set of such segments. We identify $b(\Delta)=\rho \nu_\rho^{a}$ and $e(\Delta)=\rho \nu_\rho^{b} $ with the corresponding integers $a$ and $b$. For $\Delta_1,\Delta_2 \in \mathrm{Seg}$, the condition $\Delta_1 \prec\Delta_2$ translates to
\begin{align} \label{eq:Delta1 precede Delta2}
b(\Delta_1) +1 \leq b(\Delta_2) \leq e(\Delta_1) +1 \leq e(\Delta_2).
\end{align}

Let $\Mult$ be the set of formal finite sums of segments in $\mathrm{Seg}$. There are two orders on segments in a multisegment: the \emph{left} and \emph{right aligned} orders $\le_b$ and $\le_e$. For segments $\Delta,\Delta'$, we have $\Delta \le_b \Delta'$ if $b(\Delta) \le b(\Delta')$, and either $b(\Delta)<b(\Delta')$ or $e(\Delta) \le e(\Delta')$. We have $\Delta \le_e \Delta'$ if $e(\Delta) \le e(\Delta')$, and either $e(\Delta)<e(\Delta')$ or $b(\Delta) \le b(\Delta')$. 

Note that if the segments $\Delta_i$ of a multisegment $\m$ are ordered decreasingly for either the left or right aligned order, then the mutlisegment $\m$ is ordered. Thus any multisegment can be ordered.

The Zelevinsky Classification gives a bijection
\begin{align*}
    \Mult & \to \Irr^\mathbb{Z}\\
    \underset{\text{ordered}}{\m=\Delta_1 + \cdots + \Delta_N} & \mapsto \Z(\m):=\soc(\Z(\Delta_1)\times \Z(\Delta_2)\times \cdots \times \Z(\Delta_N)).
\end{align*}

Let $\mathcal{R}$ be the Grothendieck ring of the category $\mathcal{C}^{\mathbb{Z}}$, it is freely generated, as a commutative ring, by the images of the $\Z(\Delta)$, for all $\Delta\in \mathrm{Seg}$.

\section{Quantum Affine Schur-Weyl duality} \label{sec:quantum affine Schur Weyl duality}

The correspondence between finite dimensional representations of quantum affine algebras of type $A$ and finite length representations of $p$-adic general linear groups can be formalized using quantum affine Schur-Weyl duality. 

\subsection{The functor \texorpdfstring{$\FF$}{F}}
		
The quantum Schur-Weyl duality was introduced by Jimbo in \cite{Jim86}. Let $V=\CC^k$ be the natural representation of the quantum algebra $U_q(\mathfrak{sl}_k)$. Then for all $N\geq 1$, there exists a left action of the Hecke algebra $H_N(q^2)$ on $V^{\otimes N}$ which commutes with the action of $U_q(\mathfrak{sl}_k)$. For $M$ a right $H_N(q^2)$-module, one can equip $M\otimes_{H_N(q^2)}V^{\otimes N}$ with a left $U_q(\mathfrak{sl}_k)$-module structure. Moreover, if $N < k$, the functor
\begin{equation}\label{qSWd}
M \to M\otimes_{H_N(q^2)}V^{\otimes N}
\end{equation}
is an equivalence of categories between the category of finite-dimensional $H_N(q^2)$-modules and the category of finite-dimensional $U_q(\mathfrak{sl}_k)$-modules of level $N$.

In \cite{CP96b} Chari-Pressley extended this result to the affine case. They showed that if $M$ is a right module for the affine Hecke algebra $\hat{H}_N(q^2)$, then the image of \eqref{qSWd} can be endowed with a left $U_q(\widehat{\mathfrak{sl}}_k)$-module structure. Thus defining a functor between the category of finite-dimensional modules for the affine Hecke algebra $\hat{H}_N(q^2)$ and the category $\mathscr{C}^N_{\mathfrak{sl}_k}$ of finite-dimensional $U_q(\widehat{\mathfrak{sl}}_k)$-modules of level $N$. If $k>N$, this functor is again an equivalence of categories.

On the other hand, by a result of Heiermann \cite{H11} (see also \cite{Bor76} \cite{Cas80}, \cite{BK98}, \cite{SS12} or the survey paper \cite{Gur21b}), there exists an equivalence of categories between smooth, finite-length representations of $GL_N(F)$ in a simple Bernstein block, and the category of finite-dimensional right $\hat{H}_N(q^2)$-modules. 

Thus, for $k>N$, after reduction to cuspidal lines in the $p$-adic case, and restriction of the quantum parameters for representations of the quantum affine algebra, we obtain equivalences of categories: 
\begin{equation}\label{eq:functor FkN}
\FF_{k,N}: \mathcal{C}^\ZZ(G_N) \xrightarrow{\sim} \mathscr{C}^{N,\mathbb{Z}}_{\mathfrak{sl}_k}.
\end{equation}
By summing over $N$, we get a functor $\FF$: 
\begin{equation}
\FF: \mathcal{C}^\ZZ \to \bigoplus_{N\geq 0} \mathscr{C}^{N,\ZZ}_{\mathfrak{sl}_k} \quad =: \mathscr{C}_k.
\end{equation}

%
%
%



We recall some properties of the functor $\FF$.

\begin{proposition}
\begin{enumerate}
	\item The functor $\FF$ is a tensor functor between the tensor category $\mathcal{C}^\ZZ$ and the tensor category $\mathscr{C}_k$. That is, for any pair of representations $\pi_1, \pi_2$,
	\[\FF(\pi_1 \times \pi_2) \cong \FF(\pi_1) \otimes \FF(\pi_2).\]
	\item The functor $\FF$ is exact.
	\item The functor $\FF$ send simple modules to simple modules. 
More precisely, let $\m\in \mathrm{Mult}$ a multisegment whose sum of lengths is $< k$, then $\FF_{k,N}(\Z(\m)) \cong L(M_\m)$, where the factors in the monomial are obtained from the segments of $\m$ by the following correspondence:
\begin{equation}\label{correspondance segments Ys}
\begin{array}{ccc}
 {[}a, b] & \mapsto & Y_{b-a+1, -a-b},\\
{[}\frac{1-i-p}{2}, \frac{i-p-1}{2} ] & \mapsfrom &  Y_{i,p}
\end{array}
\end{equation}

	 In particular, for $a\in \ZZ$, the image of the segment $[a,a]$ is a fundamental representation
	\[\FF(\Z([a,a])) \cong L(Y_{1,-2a})  .\]
\end{enumerate}
\end{proposition}

\subsection{Quotient and localization of category}\label{sect_quoandloc}

In this section, we parallel the results of \cite{KKK18} on the generalized quantum affine Schur-Weyl duality between representations of quantum affine algebras and those of \emph{quiver Hecke algebras} (also called \emph{Khovanov-Lauda-Rouquier}, or KLR, algebras). More specifically, a category is defined through quotient and localization of a category of finite-dimensional quiver Hecke algebra modules, for which the generalized quantum affine Schur-Weyl duality functor factors through. The resulting functor was subsequently shown by Fujita in \cite{Fuj22} to be an equivalence of categories, in all $ADE$ types.
As the KKK generalized quantum affine Schur-Weyl duality is a generalization of the type $A$ quantum affine Schur-Weyl duality, the results of \cite{KKK18} translate vertbatim to our context. 

For more details on the categorical tools used in this section, we refer the reader to \cite{KS06} (tensor categories, and localization of categories) and \cite[§ 4.3]{Pop73} (quotient categories by Serre subcategories).

\begin{proposition}[{\cite[Proposition 4.3.1, Theorem 4.3.3]{KKK18}}]\label{prop_F_segment}
Let $[a,b]$ be a segment of length $\ell=b-a+1$. Then we have
\begin{equation}
\FF(\Z([a,b])) \cong \left\lbrace \begin{array}{ll}
		L(Y_{\ell,-a-b})& \text{if } 1\leq \ell \leq k-1,\\
		\CC & \text{if } \ell=0 \text{ or } \ell=k,\\
		0 & \text{if } \ell > k.
\end{array}\right.
\end{equation}
More generally, if $\m =\Delta_1 + \cdots + \Delta_N$ is a multisegment and $\ell_1,\ldots,\ell_N$ are the corresponding lengths of the segments, then
\begin{itemize}
	\item if $\ell_i>k$ for at least one $1\leq i \leq N$, then $\FF(\Z(\m))\cong 0$,
	\item if $\ell_i \leq k$ for all $1\leq i \leq N$, then $\FF(\Z(\m))$ is a simple module.
\end{itemize}
\end{proposition}

Following \cite[Section 4]{KKK18}, we introduce a quotient of the category $\mathcal{C}^\ZZ$ in order to kill all simple modules which are send to 0 by the functor $\FF$.

Let $\mathcal{S}_k$ be the smallest Serre subcategory of $\mathcal{C}^\ZZ$ such that
\begin{enumerate}
	\item $\mathcal{S}_k$ contains $\Z([a,a+k])$, for all $a\in \ZZ$,
	\item for all $\pi_1 \in \mathcal{S}_k$ and $\pi_2 \in \mathcal{C}^\ZZ$, then $\pi_1 \times \pi_2$ and $\pi_2 \times \pi_1$ are in $\mathcal{S}_k$.
\end{enumerate}

Consider the quotient category $\mathcal{C}^\ZZ/\mathcal{S}_k$, and let $\mathcal{Q}:\mathcal{C}^\ZZ \to \mathcal{C}^\ZZ/\mathcal{S}_k$ the canonical functor. The category $\mathcal{C}^\ZZ/\mathcal{S}_k$ is an abelian tensor category. Since the functor $\FF$ send $\mathcal{S}_k$ to 0, it factors through $\mathcal{Q}$ and we obtain an exact tensor functor:
\begin{equation}
\FF': \mathcal{C}^\ZZ/\mathcal{S}_k  \to \mathscr{C}_k.
\end{equation}

We have a description of the simple objects in $\mathcal{C}^\ZZ/\mathcal{S}_k$.
\begin{proposition}[{\cite[Proposition 4.4.1]{KKK18}}]
There is a bijection:
\begin{align*}
\left\lbrace \begin{array}{c}
\text{simple objects} \\
\text{in } \mathcal{C}^\ZZ/\mathcal{S}_k
\end{array}\right\rbrace &\leftrightarrow \left\lbrace \begin{array}{c}
\text{multisegments } \Delta_1 + \cdots + \Delta_N  \\
\text{with } \ell(\Delta_i) \leq k, \forall 1 \leq i \leq N 
\end{array}\right\rbrace \\
\mathcal{Q}(\Z(\Delta_1 + \cdots + \Delta_N)) &\leftmapsto \Delta_1 + \cdots + \Delta_N.
\end{align*}
\end{proposition}

Note from Proposition~\ref{prop_F_segment}, that for all $a\in \ZZ$, the images of the $\FF'(\Z([a,a+k-1]))$ are isomorphic to the trivial representation. Thus following \cite{KKK18} we can localize $\mathcal{C}^\ZZ/\mathcal{S}_k$ to reflect this fact.

For all $a\in \ZZ$, let $Z_a=\Z([a,a+k-1])$. Then the $(Z_a)_{a\in\ZZ}$ form a commuting family of central objects in the quotient category $\mathcal{C}^\ZZ/\mathcal{S}_k$, in the sense of \cite[Appendix A]{KKK18}. 
Let $\mathcal{C}_k'$ be the localization $(\mathcal{C}^\ZZ/\mathcal{S}_k)[Z_a^{-1}\mid a\in \ZZ]$ by this commuting family. Let $\mathcal{C}_k$ denote the tensor category $(\mathcal{C}^\ZZ/\mathcal{S}_k)[Z_a\simeq 1\mid a\in \ZZ]$, as in \cite[A.7]{KKK18}. There is an exact monoidal functor $\Omega$:
\[\Omega:\mathcal{C}^\ZZ/\mathcal{S}_k \longrightarrow  \mathcal{C}_k. \]

Let us denote by $\overline{\mathcal{Q}}$ the composition of the functors $\mathcal{Q}$ and $\Omega$:
\[\overline{\mathcal{Q}}:\mathcal{C}^\ZZ \xrightarrow{Q}\mathcal{C}^\ZZ/\mathcal{S}_k \xrightarrow{\Omega}  \mathcal{C}_k. \]

Recall some properties of the functor $\overline{\mathcal{Q}}$:

\begin{proposition}[{\cite[Theorem B.1.1.-B.1.2., Lemma A.7.1.-A.7.2.]{KKK18}}]\label{proposition_properties_Qbar}
The functor $\overline{\mathcal{Q}}$ satisfies the following:
\begin{enumerate}
	\item $\overline{\mathcal{Q}}$ is exact,
	\item every exact sequence in $\mathcal{C}_k$ is isomorphic under $\overline{\mathcal{Q}}$ to the image of an exact sequence in $\mathcal{C}^\ZZ$,
	\item if $\pi\in \mathcal{S}_k$, then $\overline{\mathcal{Q}}(\pi)\simeq 0$,\label{Qbar_0}
	\item for all $a\in \mathbb{Z}$, $\overline{\mathcal{Q}}(\Z_a)\simeq {\bf 1}$, the trivial representation,\label{Qbar_1}
	\item simple objects in $\mathcal{C}_k$ are isomorphic to the images under $\overline{\mathcal{Q}}$ of simple objects in $\mathcal{C}^\ZZ$ which are not in $\mathcal{S}_k$. \label{Qbar_simples} 
\end{enumerate} 
\end{proposition}

As $k\in\ZZ_{>0}$ is fixed, we will denote by the images under $\overline{\mathcal{Q}}$ by $\overline{\pi}$, for $\pi$ a representation in $\mathcal{C}^\ZZ$, or $\overline{\Z}(\m) = \overline{\Z(\m)}$, for $\m$ a multisegment.

More precisely, we have 

\begin{corollary}\label{corollary image of Qbar}
    Let $\m=\Delta_1 + \cdots + \Delta_N\in \Mult$, and $\ell_1,\ldots,\ell_N$ the corresponding lengths of the segments. Then the image of the simple module $\Z(\m)$ in $\mathcal{C}_k$ is isomorphic to 
    \[
    \overline{\mathcal{Q}}(\Z(\m)) \cong \left\lbrace \begin{array}{cc}
        0 & \text{if } \exists i \text{ such that } \ell_i>k,\\
        \overline{\Z}(\overline{\m}^k) & \text{otherwise}
    \end{array}\right.,
    \]
where $\overline{\m}^k=\{\Delta_i\in\m \mid \ell(\Delta_i) < k \}$.
\end{corollary}

\begin{proof}
Recall that $\Z(\m)$ is defined as the socle of the representation $\Z(\Delta_1)\times \cdots \times \Z(\Delta_N) $ (assuming the multisegment $\m$ is ordered).

    If there exists $i$ such that $\ell_i>k$, then by Proposition \ref{proposition_properties_Qbar} \eqref{Qbar_0},  $\overline{\mathcal{Q}}(\Z(\Delta_1)\times \cdots \times \Z(\Delta_N)) \cong 0$ and thus $\overline{\mathcal{Q}}(\Z(\m))\cong 0$.

Otherwise, write $\overline{\m}^k= \Delta_{i_1} + \Delta_{i_2} + \cdots +\Delta_{i_s} $, then by Proposition \ref{proposition_properties_Qbar} \eqref{Qbar_1}, 
\[
\overline{\mathcal{Q}}(\Z(\Delta_1)\times \Z(\Delta_2)\times \cdots \times \Z(\Delta_N)) \cong \overline{\Z}(\Delta_{i_1})\times \overline{\Z}(\Delta_{i_2})\times \cdots \times \overline{\Z}(\Delta_{i_s}).
\]
Hence, $\overline{\mathcal{Q}}(\Z(\m)) \cong \overline{\Z}(\overline{\m}^k)$ as they are both the socle of $\overline{\Z}(\Delta_{i_1})\times \overline{\Z}(\Delta_{i_2})\times \cdots \times \overline{\Z}(\Delta_{i_s})$.
\end{proof}

In particular,
\begin{equation}
\begin{array}{rcll}
\left\lbrace \begin{array}{c}
\text{simple objects} \\
\text{in } \mathcal{C}_k
\end{array}\right\rbrace & \leftrightarrow & \left\lbrace  \begin{array}{c}
\text{multisegments } \Delta_1 + \cdots + \Delta_N  \\
\text{with } \ell(\Delta_i) < k, \forall 1 \leq i \leq N 
\end{array}\right\rbrace & =:\Mult_k \\
\overline{\mathcal{Q}}(\Z(\Delta_1 + \cdots + \Delta_N)) &\leftmapsto& \Delta_1 + \cdots + \Delta_N.
\end{array}
\end{equation}

We have the following.

\begin{theorem}[{\cite[Theorem 4.6.6, Proposition A.7.3]{KKK18}}]
The functor $\FF':\mathcal{C}^\ZZ/\mathcal{S}_k \to \mathscr{C}_k$ factors through $\mathcal{C}_k$, leading to an exact functor:
\begin{equation}
\tilde{\FF} : \mathcal{C}_k \to \mathscr{C}_k.
\end{equation}
\end{theorem}

Based on similar results between representations of quantum affine algebras and representations of quiver Hecke algebras (see \cite{Fuj22}), we formulate the following: 

\begin{conjecture}\label{conj: equivalence of categories}
The functor $\tilde{\FF}$ is an equivalence of categories.
\end{conjecture}

\subsection{Isomorphism of Grothendieck rings}\label{sect:Grothendieck ring approach}


From now on, let us drop the superscript $\ZZ$ for clarity of notations. 

For all $k\geq 1$, let $\mathcal{R}_k$ be the Grothendieck ring of the category $\mathcal{C}_k$. The following is result provides some evidence toward Conjecture \ref{conj: equivalence of categories}.

\begin{prop}\label{prop: iso Rk's}
The exact functor $\tilde{\FF}$ induces a ring isomorphism between the Grothendieck rings:
\[ \mathcal{R}_k \cong \mathscr{R}_k \quad (=K_0(\rep(U_q(\hat{\mathfrak{sl}}_k))).
\]
\end{prop}

\begin{proof}
From the previous results, the functor $\tilde{\FF}$ is monoidal and exact. Let us prove that $\tilde{\FF}$ induces a bijection between the simple objects in $\mathcal{C}_k$ and $\mathscr{C}_k$.

Let $S$ be a simple object in $\mathcal{C}_k$, then there exists a multisegment $\m\in \Mult_k$, such that $S\cong \overline{\Z}(\m)$, and by Proposition \ref{prop_F_segment}, $\tilde{\FF}(S)$ is simple. 
Conversely, let $L$ be simple object in $\mathscr{C}_k$, then it can be obtained as the socle of a tensor product of fundamental representations:
\[
L=\soc(L(Y_{i_1,r_1})\otimes L(Y_{i_2,r_2})\otimes \cdots \otimes L(Y_{i_s,r_s})).
\]
Using correspondence \eqref{correspondance segments Ys}, define, for all $1\leq j\leq s$, $a_j=\frac{1-i_j-r_j}{2}$ and $b_j=\frac{i_j-r_j-1}{2}$, then
\[
\tilde{F}(\overline{\Z}([a_1,b_1])\times \overline{\Z}([a_2,b_2])\times \cdots \overline{\Z}([a_s,b_s]) ) \cong L(Y_{i_1,r_1})\otimes L(Y_{i_2,r_2})\otimes \cdots \otimes L(Y_{i_s,r_s}). 
\]
Then $\tilde{\FF}(\overline{\Z}([a_1,b_1] + [a_2,b_2]+ \cdots + [a_s,b_s]))$ is isomorphic to the socle of $L(Y_{i_1,r_1})\otimes L(Y_{i_2,r_2})\otimes \cdots \otimes L(Y_{i_s,r_s})$ and hence to $L$. 
Thus, $\tilde{\FF}$ induces a ring isomorphism from $\mathcal{R}_k$ to $\mathscr{R}_k$.


\end{proof}

\begin{remark}
    As a consequence of Corollary \ref{corollary image of Qbar}, for $\m = \Delta_1 + \Delta_2 + \cdots + \Delta_N$ a multisegment, its image in $\mathcal{R}_k$ is
\begin{equation}
[\Z(\m)]_k = \left\lbrace\begin{array}{ccl}
&0 & \text{ if } \m \notin\mathrm{Mult}_{k+1}, \\
&[\overline{\Z}(\overline{\m}^k)] & \text{ otherwise}
\end{array}\right..
\end{equation}
\end{remark}

\section{Quantum affine and \texorpdfstring{$k$}{k}-Zelevinsky classifications} \label{sec: QA and k-Z classification}


In this section, we establish analogs of the Zelevinsky classification \cite{Zel}, in terms of representations of the type $A$ quantum affine algebra, and in the category $\mathcal{C}_k$.


Recall from Section \ref{sec:quantum affine Schur Weyl duality}, under the equivalence of categories, multisegments and dominant monomials are identified via the following correspondence between segments and fundamental monomials:

\begin{equation}\label{eq:multisegtomonom}
\begin{array}{ccc}
\mathrm{Seg}_k = \{ [a,b] \mid b-a+1 <k \} & \simeq & \mathcal{P}^+ = \{ Y_{i,p} \mid (i,p)\in \hat{I} \},\\
 {[}a, b] & \mapsto & Y_{b-a+1, -a-b},\\
{[}\frac{1-i-p}{2}, \frac{i-p-1}{2} ] & \mapsfrom &  Y_{i,p}
\end{array}
\end{equation}
We denote this correspondence by ${\bf m} \mapsto M_{\bf m}$ and $M \mapsto {\bf m}_M$ accordingly. 

\begin{remark}
The correspondence between dominant monomials and multisegments we use here is different from the one used in \cite{CDFL}. The one used in \cite{CDFL} is given by $[a, b] \mapsto Y_{b-a+1, a+b-1}$. Accordingly, the correspondence between dominant monomials and semistandard Young tableaux used in this paper is also different from the one used in \cite{CDFL}. 
\end{remark}

\subsection{Quantum affine Zelevinsky Classification}

\begin{lemma}\label{lemma_one_neg}
Let $(i,p)\in \hat{I}$, the only monomials appearing in $\chi_q(L(Y_{i,p}))$ with exactly one negative power of $Y$ have the form
\begin{equation}\label{eq_one_neg_mon}
Y_{q,p+i-q}\left( Y_{q+r,p+i+r-q}\right)^{-1}Y_{i+r,p+r},
\end{equation}
where $\left\lbrace\begin{array}{l}
0\leq q \leq i-1,\\
1\leq r\leq k-i
\end{array}\right.$.
\end{lemma}

\begin{proof}
From the work of Nakajima, explicit formulas for the $q$-characters of fundamental representations are know, in type $A$. Let $(i,p)\in \hat{I}$, then by \cite[Proposition 4.6]{N03}:
\begin{equation}\label{eq_chiqLY}
\chi_q(L(Y_{i,p})) = \sum_{1\leq j_1 < \cdots  < j_i \leq k} \prod_{m =1}^i Y_{j_m, p+i + j_m-2m}\left( Y_{j_m-1, p+i + j_m-2m +1} \right)^{-1}.
\end{equation}

Fix $1\leq j_1 < \cdots  < j_i \leq k$, and let $m$ be the corresponding monomial of $\chi_q(L(Y_{i,p}))$ in \eqref{eq_chiqLY}. We note that if $j_{m+1}=j_m+1$, then the negative powered $Y$-variable in the $(m+1)$th factor of $m$ will cancel out with the positive powered $Y$-variable in the $m$th factor of $m$. 
Hence, the monomial $m$ has exactly one negative powered $Y$-variable if and only if the tuple $(j_1, j_2, \cdots, j_i)$ is formed of two disconnected intervals of integers, thus if and only if it has the form
\[(1,2,\ldots,q, q+r+1, q+r+2,  \ldots, r+i ),\]
with $0\leq q \leq i-1$, $r\leq 1$ and $r+i\leq k$.
In that case, the monomial is the one from \eqref{eq_one_neg_mon}.
\end{proof}

\begin{proposition}\label{prop_red_ineq_j_s}
Let $(i,p),(j,s)\in \hat{I}$ such that $p\leq s$. The tensor product $L(Y_{j,s})\otimes L(Y_{i,p})$ is reducible if and only if 
\begin{align}
\label{ineq_j+s} p+i+2 \leq & j+s \leq 2k +p-i,\\
\label{ineq_j-s} -p-i \leq & j-s \leq i-p-2.
\end{align}
Moreover in that case, the tensor product is of length 2 and we have the following short exact sequence
\begin{equation}\label{eq_SES}
0 \to L(Y_{q,p+i-q}Y_{i+r,p+r}) \to L(Y_{j,s})\otimes L(Y_{i,p}) \to L(Y_{j,s}Y_{i,p}) \to 0,
\end{equation}
where $\left\lbrace\begin{array}{l}
q= \frac{1}{2}(i+j+p-s),\\
r= \frac{1}{2}(-i+j-p+s)
\end{array}\right. $.
\end{proposition}

\begin{proof}
Let us consider the dominant monomials appearing in the product of $q$-characters $\chi_q(L(Y_{j,s}))\chi_q(L(Y_{i,p}))$. One of these dominant monomials is naturally the product of the highest weights $Y_{j,s}Y_{i,p}$. As $p\leq s$, the only other possible dominant monomials are of the form $Y_{j,s}M$, with $M$ a monomial of $\chi_q(L(Y_{i,p}))$ with exactly one negative powered $Y$-variable.
From Lemma~\ref{lemma_one_neg}, we know that these $M$ have the form 
\begin{equation}\label{eq_ineq_q_r}
M=Y_{q,p+i-q}\left( Y_{q+r,p+i+r-q}\right)^{-1}Y_{i+r,p+r}, \quad \text{with} \quad \left\lbrace \begin{array}{l}
0\leq  q \leq i-1,\\
1\leq r\leq k-i
\end{array}\right. .
\end{equation}

Thus $Y_{j,s}M$ is dominant if and only if we have $Y_{j,s}=Y_{q+r,p+i+r-q}$, or equivalently
\begin{align*}
j&=q+r ,\\
s&=p+i+r-q.
\end{align*}

Suppose now that the tensor product $L(Y_{j,s})\otimes L(Y_{i,p})$ is reducible, then necessarily one of the $Y_{j,s}M$ is dominant. Thus $(j,s)$ satisfy the inequalities of \eqref{eq_ineq_q_r}, which are equivalent to the inequalities \eqref{ineq_j+s} and \eqref{ineq_j-s}.

Conversely, suppose now that the inequalities \eqref{ineq_j+s}, \eqref{ineq_j-s} are satisfied. In particular, one has 
\[s-p \geq |i-j| +2,\]
and thus $L(Y_{j,s}Y_{i,p})$ is a \emph{snake module}, and using \cite[Theorem 6.1]{MY12}, we know $Y_{j,s}Y_{i,p}$ is its unique dominant monomial. 

Let us write $j=q+r$ and $s=p+i+r-q$, then $(q,r)$ satisfy the inequalities in \eqref{eq_ineq_q_r}. Necessarily, the corresponding monomial $Y_{j,s}M$ is dominant; it is also the only one of this form. Thus we have the following relation between the $q$-characters
\begin{equation}\label{eq_prod_chiq_LL}
\chi_q(L(Y_{j,s}))\chi_q(L(Y_{i,p})) = \chi_q(L(Y_{j,s}Y_{i,p})) + \chi_q(L(Y_{q,p+i-q}Y_{i+r,p+r})).
\end{equation}

In particular, $\chi_q(L(Y_{j,s}Y_{i,p})) \neq \chi_q(L(Y_{j,s}))\chi_q(L(Y_{i,p}))$ and the tensor product $L(Y_{j,s})\otimes L(Y_{i,p})$ is reducible.

Moreover, from \cite[Theorem 9.1]{K02} (see also \cite[Corollary 5.1]{C02}), as $s\geq p$, $L(Y_{j,s}\otimes Y_{i,p})$ is a cyclic module, generated by the tensor product of the highest weight vectors of its factors. Hence it has a unique maximal submodule and a unique simple quotient. From the highest weight, its unique simple quotient is isomorphic to $L(Y_{j,s}Y_{i,p})$. Thus by \eqref{eq_prod_chiq_LL}, its unique maximal proper submodule has the $q$-character of $L(Y_{q,p+i-q}Y_{i+r,p+r})$. Since $L(Y_{q,p+i-q}Y_{i+r,p+r})$ is simple, it is isomorphic to the unique maximal proper submodule of $L(Y_{j,s}\otimes Y_{i,p})$. In conclusion we have the short exact sequence \eqref{eq_SES}.
\end{proof}

\begin{remark}\label{rem_intro prec k}
In the context of $p$-adic representations of $\mathrm{GL}_N$, as a consequence of the Zelevinsky classification, if $\Delta,\Delta'$ are unlinked, then $\Z(\Delta)\times\Z(\Delta')$ is irreducible. Moreover, if $\Delta\prec\Delta'$, then $\Z(\Delta)\times\Z(\Delta')$ is of length 2, and we have the following short exact sequence (see \cite[Lemma A.9]{LM16})
\begin{equation*}
0 \to \Z(\Delta \cup \Delta')\times \Z(\Delta \cap \Delta')\to \Z(\Delta)\times\Z(\Delta') \to \Z(\Delta+\Delta')\to 0.
\end{equation*}
If through the equivalence \eqref{eq:multisegtomonom}, $\Delta \mapsto Y_{j,s}$ and $\Delta' \mapsto Y_{i,p}$, this is the same short exact sequence as \eqref{eq_SES}. Moreover, we note that the conditions \eqref{ineq_j+s} and \eqref{ineq_j-s} translate  to:
\begin{equation*}
    \Delta \prec \Delta' \quad \text{and} \quad e(\Delta') - b(\Delta) \leq k-1.
\end{equation*}

\end{remark}

For representations of the quantum affine algebra $U_q(\widehat{\mathfrak{sl}}_k)$, the irreducible representation $L(M)$ of highest weight $M\in \mathcal{P}$ is usually defined as the quotient of the Verma module of the corresponding weight. However, it can also be characterized as the socle of a tensor product of fundamental representations. 

\begin{proposition}\cite{AK97}
Let $m=Y_{i_1,p_1}Y_{i_2,p_2}\cdots Y_{i_N,p_N}$, such that for all $1\leq j\leq N$, $(i_j,p_j)\in\hat{I}$ and $p_1\leq p_2 \leq \cdots \leq p_N$, then $L(m)$ is isomorphic to the unique irreducible submodule (hence the socle) of the \emph{standard module}
\begin{equation}
M(m):= L(Y_{i_1,p_1})\otimes L(Y_{i_2,p_2})\otimes \cdots \otimes L(Y_{i_N,p_N}).
\end{equation}
\end{proposition}



Recall that in the Zelevinsky classification in Section~\ref{sec:quantum affine algebras and p-adic groups}, the irreducible representation $\Z(\m)$ corresponding to a multisegment $\m= \Delta_1 + \cdots + \Delta_N$ is defined as the socle of the parabolic induction $\Z(\Delta_1)\times \cdots \times \Z(\Delta_N)$, where the $\Delta_j$ are ordered such that for all $j<\ell$, $\Delta_j$ does not precede $\Delta_\ell$. Through the correspondence \eqref{eq:multisegtomonom}, the same condition is applied here. Indeed, for all $j<\ell$, the segment $\left[\frac{1-i_j-p_j}{2}, \frac{i_j-p_j-1}{2}\right]$ does not precede the segment $\left[\frac{1-i_\ell-p_\ell}{2}, \frac{i_\ell-p_\ell-1}{2}\right]$.



\subsection{Zelevinsky classification in the category \texorpdfstring{$\mathcal{C}_k$}{Ck}} \label{subsec:Zelevinsky classification in Ck}

Inspired by Remark~\ref{rem_intro prec k}, we introduce the following.

\begin{definition}
If $\Delta,\Delta'$ are two segments of length $<k$. We say that $\Delta$ \emph{$k$-precedes} $\Delta'$ if $\Delta\prec\Delta'$ and $e(\Delta') - b(\Delta) <k$. We denote it by $\Delta\prec_k\Delta'$. If either $\Delta\prec_k\Delta'$ or $\Delta'\prec_k\Delta$ we say that $\Delta$ and $\Delta'$ are \emph{$k$-linked}.
\end{definition}

\begin{lemma}\label{lem_Delta_Delta'_k_unlinked}
For $\Delta, \Delta'$ segments in $\Mult_k$, the product $\overline{\Z}(\Delta)\times \overline{\Z}(\Delta')$ is simple in $\mathcal{C}_k$ if and only if $\Delta$ and $\Delta'$ are $k$-unlinked.
\end{lemma}

\begin{proof}
Recall from Remark \ref{rem_intro prec k} that by \cite[Lemma A.9]{LM16}, if $\Delta\prec\Delta'$, then $\Z(\Delta)\times\Z(\Delta')$ is of length 2, and we have the following short exact sequence
\begin{equation}\label{eq_cup_cap}
0 \to \Z(\Delta \cup \Delta')\times \Z(\Delta \cap \Delta')\to \Z(\Delta)\times\Z(\Delta') \to \Z(\Delta+\Delta')\to 0.
\end{equation}

If moreover $\Delta\prec_k\Delta'$, then $\overline{\mathcal{Q}}(\Z(\Delta \cup \Delta')\times \Z(\Delta \cap \Delta'))$ is non-zero, thus $\overline{\Z}(\Delta \cup \Delta')\times \overline{\Z}(\Delta \cap \Delta')$ is a proper subobject of $\overline{\Z}(\Delta)\times \overline{\Z}(\Delta')$, which is not simple.

If $e(\Delta')-b(\Delta)\geq k$, then the length of $\Delta \cup \Delta'$ is $\geq k+1$ and thus $\overline{\Z}(\Delta \cup \Delta')=0$. The image of Equation \eqref{eq_cup_cap} in $\mathcal{C}_k$ gives that $\overline{\Z}(\Delta)\times\overline{\Z}(\Delta')$ is simple and isomorphic to $\overline{\Z}(\Delta+\Delta')$.

If $\Delta$ and $\Delta'$ are unlinked, then $\Z(\Delta)\times\Z(\Delta')$ is irreducible and thus $\overline{\Z}(\Delta)\times \overline{\Z}(\Delta')$ is simple.
\end{proof}

For $\m=\Delta_1 + \cdots + \Delta_N\in \Mult_k$ an ordered multisegment, let us write
\[\overline{\zeta}(\m):=\overline{\Z}(\Delta_1)\times \cdots \times \overline{\Z}(\Delta_N).\]

The following is the key element of Zelevinsky's classification,  it generalizes Lemma~\ref{lem_Delta_Delta'_k_unlinked}.

\begin{proposition}
For $\m = \Delta_1 + \Delta_2 + \cdots + \Delta_N\in\Mult_k$, $\overline{\zeta}(\m)$ is simple if and only if the segments $\Delta_i$ are pairwise $k$-unlinked.
\end{proposition}

We use the following intermediate result\footnote{In \cite{Zel}, the proof is given in terms of Langlands classification but it is actually valid for both classifications.}.

\begin{lemma}[{\cite[Section 7]{Zel}}]\label{lemma_JH_zeta}
The terms in the Jordan-Hölder series of $\Z(\Delta_1)\times \cdots \times \Z(\Delta_N)$ are of the form $\Z(\n)$, where $\n$ is obtained from $\m$ by a sequence of operations where a pair of linked segments $\{\Delta,\Delta'\}$ are replaced by $\{\Delta\cap \Delta',\Delta\cup\Delta'\}$. 
\end{lemma}

\begin{proof}
If the segments $\Delta_i$ are pairwise $k$-unlinked, for all pair of possibly linked segments $\Delta_i,\Delta_j$, the union $\Delta_i\cup \Delta_j$ is of length $\geq k+1$. Thus by Lemma~\ref{lemma_JH_zeta}, all terms in the Jordan-Hölder series of $\Z(\Delta_1)\times \cdots \times \Z(\Delta_N)$ are sent to 0 by $\overline{\mathcal{Q}}$ except $\Z(\m)$. Hence $\overline{\zeta}(\m)$ is simple and isomorphic to $\overline{\Z}(\m)$.

Conversely, if $\overline{\zeta}(\m)$ is simple, then all terms in the Jordan-Hölder series of $\Z(\Delta_1)\times \cdots \times \Z(\Delta_N)$ are sent to 0 by $\overline{\mathcal{Q}}$. We deduce that the segments of $\m$ are pairwise $k$-unlinked.
\end{proof}

\begin{corollary}[$k$-Zelevinsky classification]
For $\m = \Delta_1 + \Delta_2 + \cdots + \Delta_N\in\Mult_k$ such that $\Delta_i\nprec_k\Delta_j$ for all $i<j$, 
\[\overline{\Z}(\m) = \soc(\overline{\Z}(\Delta_1)\times \overline{\Z}(\Delta_2)\times \cdots \times \overline{\Z}(\Delta_N)) \quad \text{ in } \mathcal{C}_k.\]
\end{corollary}

The following is a consequence of the $k$-Zelevinsky classification.
\begin{lemma}\label{lemma_m_ncpreck_n}
Write $\m=\sum_i\Delta_i$ and $\n=\sum_j\Delta'_j$. If for all $i,j$, $\Delta_i\nprec_k \Delta'_j$, then $\overline{\Z}(\m + \n)= \soc(\overline{\Z}(\m) \times \overline{\Z}(\n)) = \coso(\overline{\Z}(\n) \times \overline{\Z}(\m)) $ in $\mathcal{C}_k$.
\end{lemma}

\section{A criterion for the simplicity of tensor products} \label{sec:a criterion for simplicity of tensor products}

\subsection{Condition for irreducibility} \label{subsec:LI,RI}

Following \cite{LM16}, for $\pi = \overline{\rm Z}(\m)$, $\sigma = \overline{\rm Z}({\bf n})$, irreducible objects in $\mathcal{C}_k$, denote by ${\rm LI}_k(\pi, \sigma)$ the condition $\overline{\rm Z}(\m +\n) = {\rm soc}(\pi \times \sigma)$ in $\mathcal{C}_k$, and denote by ${\rm RI}_k(\pi, \sigma)$ the condition $\overline{\rm Z}(\m + \n) = {\rm cos}(\pi \times \sigma)$ in $\mathcal{C}_k$. We also use the notations from \cite{LM16}. For $\pi = {\rm Z}(\m)$, $\sigma = {\rm Z}({\bf n})$, irreducible objects in $\mathcal{C}^\ZZ$, denote by ${\rm LI}(\pi, \sigma)$ the condition ${\rm Z}(\m +\n) = {\rm soc}(\pi \times \sigma)$ in $\mathcal{C}^\ZZ$, and denote by ${\rm RI}(\pi, \sigma)$ the condition ${\rm Z}(\m + \n) = {\rm cos}(\pi \times \sigma)$ in $\mathcal{C}^\ZZ$.

\begin{lemma} [{\cite[Lemma 4.2]{LM16}}]\label{lem_LI_RI}
For $\pi={\rm Z}(\m), \sigma={\rm Z}(\n) \in {\rm Irr}_k$, $\pi \times \sigma$ is irreducible if and only if ${\rm LI}_k(\pi, \sigma)$ and ${\rm RI}_k(\pi, \sigma)$. 
\end{lemma}

The following is a consequence of Proposition~\ref{proposition_properties_Qbar}.
\begin{proposition}\label{proposition_LI_implies_LIk}
For $\pi ={\rm Z}(\m)$, $\sigma = {\rm Z}({\bf n})$, irreducible representations in $\mathcal{C}^\ZZ$, such that all segments in $\m$ and $\n$ have length inferior to $k$, then we have the implication
\begin{equation}
{\rm LI}(\pi,\sigma) \Rightarrow {\rm LI}_k(\overline{\pi},\overline{\sigma}).
\end{equation}
\end{proposition}

\begin{remark}
Lemma~\ref{lemma_m_ncpreck_n} translates as the following: if $\m=\sum_i\Delta_i$ and $\n=\sum_j\Delta'_j$ satisfy $\Delta_i\nprec_k \Delta'_j$, for all $i,j$, then ${\rm LI}_k(\overline{\Z}(\m),\overline{\Z}(\n))$ holds. 
\end{remark}

\subsection{Matching functions} \label{subsec:matching functions}


Following \cite{LM16}, let $X, Y$ be finite sets and $\rightsquigarrow$ a relation between $X$ and $Y$. An injective function $f: X \to Y$ satisfying $f(x) \rightsquigarrow x$ for all $x \in X$ is called a $\rightsquigarrow$-matching function between $X$ and $Y$. An injective function $f$ from a subset of $X$ to $Y$ satisfying $f(x) \rightsquigarrow x$ for all $x$ in the domain of $f$ is called a $\rightsquigarrow$-matching between $X$ and $Y$. 
Suppose that $X$ and $Y$ are totally ordered with respect to $\le_X$ and $\le_Y$ respectively. Define a $\rightsquigarrow$-matching $f$ between $X$ and $Y$ and its domain $I$ recursively by $x \in I$ if and only if there is $y \in Y\backslash f(I \cap X_{>x})$ such that $y \rightsquigarrow x$ in which case $f(x) = \min\{ y \in Y\backslash f(I \cap X_{>x}): y \rightsquigarrow x \}$. This $\rightsquigarrow$-matching is called a best $\rightsquigarrow$-matching between $X$ and $Y$. 

For two multisegments ${\bf m}=\sum_{i=1}^N \Delta_i$, ${\bf n} = \sum_{i=1}^{N'} \Delta'_i \in {\rm Mult}$, and two sets $X, Y \subset [N] \times [N']$, define a relation $\rightsquigarrow_{{\bf m}, {\bf n}}$ (or simply $\rightsquigarrow$ if there is no confusion) between $Y$ and $X$ as follows: $(i_2, j_2)  \rightsquigarrow (i_1, j_1)$ if and only if either $i_1=i_2$, $\Delta'_{j_2} \prec \Delta'_{j_1}$ or $j_1=j_2$, $\Delta_{i_1} \prec \Delta_{i_2}$. 

For two multisegments ${\bf m}=\sum_{i=1}^N \Delta_i$, ${\bf n} = \sum_{i=1}^{N'} \Delta'_i \in {\rm Mult}_k$, and two sets $X, Y \subset [N] \times [N']$, define a relation $\rightsquigarrow_{k, {\bf m}, {\bf n}}$ (or simply $\rightsquigarrow_k$ if there is no confusion) between $Y$ and $X$ as follows: $(i_2, j_2)  \rightsquigarrow_k (i_1, j_1)$ if and only if either $i_1=i_2$, $\Delta'_{j_2} \prec_k \Delta'_{j_1}$ or $j_1=j_2$, $\Delta_{i_1} \prec_k \Delta_{i_2}$. 

\subsection{The combinatorial conditions \texorpdfstring{${\rm LC}(\m, \m')$ and ${\rm LC}_k(\m, \m')$ }{LC(m,m') and LCk(m,m') }} \label{subsec:LCmn, LCkmn}

For two multisegments ${\bf m}, {\bf m}'$, Lapid and Minguez \cite{LM16} defined 
\begin{align*}
& X_{{\bf m}, {\bf m}'} = \{ (i,j): \Delta_i \prec \Delta'_j \}, \\
& Y_{{\bf m}, {\bf m}'} = \{ (i,j): \overleftarrow{\Delta_i} \prec \Delta'_j \}.
\end{align*}
They \cite{LM16} introduced a condition ${\rm LC}({\bf m}, {\bf m}')$: there is an injective map $f: X_{{\bf m}, {\bf m}'} \to Y_{{\bf m}, {\bf m}'}$, $f(i,j)=(i',j')$, such that for any $(i,j) \in X_{{\bf m}, {\bf m}'}$, either $i=i'$, $\Delta'_{j'} \prec \Delta'_{j}$ or $j=j'$, $\Delta_{i} \prec \Delta_{i'}$.

We consider the analogs for the quotient ring of the sets used by Lapid and Minguez \cite{LM16} to establish a combinatorial criterion to determine the irreducibly of a product of representations. 
For $k \in \mathbb{Z}_{\ge 1}$ and two multisegments ${\bf m}=\sum_{i} \Delta_i$, ${\bf m}'=\sum_{j} \Delta'_j$,  
\begin{align*}
& X^{(k)}_{{\bf m}, {\bf m}'} = \{(i,j): \Delta_i \prec_k \Delta'_j \}, \\
& Y^{(k)}_{{\bf m}, {\bf m}'} = \{(i,j): \overleftarrow{\Delta_i} \prec_k \Delta'_j \}.
\end{align*}

We denote by ${\rm LC}_k({\bf m}, {\bf m}')$ the following condition: there is an injective map $f: X^{(k)}_{{\bf m}, {\bf m}'} \to Y^{(k)}_{{\bf m}, {\bf m}'}$, $f(i,j)=(i',j')$, such that for any $(i,j) \in X^{(k)}_{{\bf m}, {\bf m}'}$, either $i=i'$, $\Delta'_{j'} \prec_k \Delta'_{j}$  or $j=j'$ and $\Delta_{i} \prec_k \Delta_{i'}$.

\begin{lemma} \label{lem:LC implies LCk}
If ${\bf m}=\sum_{i} \Delta_i$ and ${\bf m}'=\sum_{j} \Delta'_j$ are multisegment, and $k\in\mathbb{Z}_{\geq 1}$, the following implications hold
\begin{align}
 {\rm LC}_{k+1}({\bf m}, {\bf m}') &\quad \Rightarrow \quad {\rm LC}_k({\bf m}, {\bf m}') ,\label{eq:LC_k+1 => LC_k}\\
{\rm LC}({\bf m}, {\bf m}') & \quad \Rightarrow \quad  {\rm LC}_k({\bf m}, {\bf m}')\label{eq:LC => LC_k} .
\end{align}
\end{lemma}

\begin{proof}
Suppose ${\rm LC}({\bf m}, {\bf m}')$ and let $f$ be the matching function between the sets $ X_{{\bf m}, {\bf m}'}$ and $Y_{{\bf m}, {\bf m}'}$. Let $\tilde{f}$ be the restriction of $f$ to the subset $X^{(k)}_{{\bf m}, {\bf m}'}$. Then we can see that the image of $\tilde{f}$ is contained in $Y^{(k)}_{{\bf m}, {\bf m}'}$, and $\tilde{f}$ defines a matching:
\[\tilde{f}: X^{(k)}_{{\bf m}, {\bf m}'} \to  Y^{(k)}_{{\bf m}, {\bf m}'}.\] 
Indeed, let $(i,j)\in X^{(k)}_{{\bf m}, {\bf m}'}$ and $\tilde{f}(i,j)=(i',j')$ be its image. If $i=i'$, then by definition of $f$, $\Delta'_{j'}\prec \Delta'_j$, and thus $e(\Delta'_{j'})\leq e(\Delta'_j)-1$. Then $e(\Delta'_{j'}) - b(\Delta_i) \leq e(\Delta'_j) -1 - b(\Delta_i) < k-1$, as $(i,j)\in X^{(k)}_{{\bf m}, {\bf m}'}$. Thus $(i,j')\in Y^{(k)}_{{\bf m}, {\bf m}'}$. Moreover, as $\overleftarrow{\Delta_i}\prec \Delta'_{j'}$, we also have $e(\Delta'_j) - b(\Delta'_{j'}) < k$. 
If $j=j'$, we show similarly that $(i',j)\in Y^{(k)}_{{\bf m}, {\bf m}'}$ and that $e(\Delta_{i'}) - b(\Delta_i) < k$.

The implication \eqref{eq:LC_k+1 => LC_k} can be proven in a similar way, as it is clear that $X^{(k)}_{{\bf m}, {\bf m}'} \subset X^{(k+1)}_{{\bf m}, {\bf m}'}$ and $Y^{(k)}_{{\bf m}, {\bf m}'} \subset Y^{(k+1)}_{{\bf m}, {\bf m}'}$.
\end{proof}

\subsection{Main result}

The following is the analog in our context of {\cite[Proposition 6.20]{LM16}}.

\begin{proposition}\label{proposition_LC_LI_cusp}
Let $\n \in {\rm Mult}_k$, and $a\in\mathbb{Z}$, then
\begin{enumerate}
\item The condition ${\rm LC}_k([a,a],\n)$ is equivalent to ${\rm LI}_k(\Z([a,a]),\Z(\n))$.

\item  The condition ${\rm RC}_k([a,a],\n)$ is equivalent to ${\rm RI}_k(\Z([a,a]),\Z(\n))$.
\end{enumerate}
\end{proposition}

We conjecture the following generalization to Proposition~\ref{proposition_LC_LI_cusp}.

\begin{conjecture}
Proposition~\ref{proposition_LC_LI_cusp} is also true if we replace $[a,a]$ by any ladder $\m$.
\end{conjecture}

Proposition \ref{proposition_LC_LI_cusp} will be proven in Section~\ref{subsec:the case that m is a segment with length 1}.

\begin{proposition}\label{prop: ladder times ladder case}
Let $\m,\n\in \mathrm{Mult}_k$ be ladders, the following conditions are equivalent:
\begin{enumerate}
	\item ${\rm LC}_k(\m, \n)$ and ${\rm LC}_k(\n, \m)$,
	\item $[\overline{\Z}(\m)\times \overline{\Z}(\n)] = [\overline{\Z}(\m + \n)] \in \mathcal{R}_k$.
\end{enumerate}
\end{proposition}

Proposition \ref{prop: ladder times ladder case} will be proven in Section~\ref{sec:proof of the case of ladder times ladder}. We can now prove the main result.

\begin{theorem} \label{thm:condition of irreducibility of tensor product}
For any simple snake $U_q(\widehat{\mathfrak{sl}_k})$-modules $L(M), L(M')$, the following conditions are equivalent:
\begin{enumerate}
\item ${\rm LC}_k(\m_M,\m_{M'})$ and ${\rm LC}_k(\m_{M'}, \m_M)$,
\item the tensor product $L(M)\otimes L(M')$ is irreducible.
\end{enumerate} 
%
%
%
%
\end{theorem}

\begin{proof}
If both $L(M)$ and $L(M')$ are snake modules, then from Proposition~\ref{prop: ladder times ladder case}, $(1)$ is equivalent to the condition $[\overline{\Z}(\m_M)\times \overline{\Z}(\m_{M'})] = [\overline{\Z}(\m_M + \m_{M'})] \in \mathcal{R}_k$. Using Proposition~\ref{prop: iso Rk's}, the latter is equivalent to the condition $[L(M)\otimes L(M')] = [L(MM')] \in \mathscr{R}_k$, which can also be translated to $\chi_q(L(M))\chi_q(L(M'))=\chi_q(L(M M'))$. Finally, this condition is equivalent to $(2)$.
\end{proof}

\begin{proposition}\label{prop:condition of irred for Yk-1}
Theorem~\ref{thm:condition of irreducibility of tensor product} is also true if $L(M')$ is a fundamental representation at an extremal node $L(Y_{1,r})$ or $L(Y_{k-1},r)$, and $L(M)$ any simple module.
\end{proposition}

\begin{proof}
If $L(M')$ is the fundamental representation $L(Y_{1,-2a})$, and $L(M)$ is any simple module, then using Proposition~\ref{proposition_LC_LI_cusp}, the condition $(1)$ of Theorem~\ref{thm:condition of irreducibility of tensor product} is equivalent to both the conditions ${\rm LI}_k(\Z([a,a]),\Z(\m_M))$ and ${\rm RI}_k(\Z([a,a]),\Z(\m_M))$ being satisfied. From Lemma~\ref{lem_LI_RI}, these conditions are equivalent to $\Z([a,a])\times\Z(\m_M)$ being irreducible in $\mathcal{C}_k$, and thus to $[\overline{\Z}([a,a])\times\overline{\Z}(\m_M)] = [\overline{\Z}([a,a] +\m_M )] \in \mathcal{R}_k$. With the same reasoning as above, this condition is equivalent to $(2)$ of Theorem~\ref{thm:condition of irreducibility of tensor product}.

The result carries to the case of the fundamental representation $L(Y_{k-1},r)$ as it is dual to a fundamental module by symmetry of the Dynkin diagram (see \cite{FM01}). Note that, as for the segments of length 1, if $\m = [a,a+k-2]$ and $\n$ is any multisegment, then $X_{\m,\n} = X_{\m,\n}^{(k)}$ and $Y_{\m,\n} = Y_{\m,\n}^{(k)}$, and the proof of Proposition~\ref{proposition_LC_LI_cusp} carries to that case.
\end{proof}



Considering the results of \cite{LM18}, we expect the following.

\begin{conjecture} \label{conj:equivalence of (2) and (3) in the case that one is snake module}
Theorem~\ref{thm:condition of irreducibility of tensor product} is true if either $L(M)$ or $L(M')$ is a snake module and the other is an arbitrary simple module.
\end{conjecture}

We give two examples to explain Theorem \ref{thm:condition of irreducibility of tensor product}.

\begin{example} \label{example:fundamental times fundamental}
Let ${\bf m} = [-4,-1]$, ${\bf m}' = [-1,2]$. Then $X_{{\bf m}, {\bf m}'} = Y_{{\bf m}, {\bf m}'} = \{(1,1)\}$, $X_{{\bf m}', {\bf m}} = Y_{{\bf m}', {\bf m}} = \emptyset$. The condition ${\rm LC}( {\bf m}, {\bf m}' )$ is not satisfied and the condition ${\rm LC}( {\bf m}', {\bf m} )$ is satisfied. Therefore the representation ${\rm Z}({\bf m}) \times {\rm Z}({{\bf m}'})$ in $\mathcal{C}$ is not irreducible.

In order for $[\overline{\Z}(\m)]$ and $[\overline{\Z}(\m')]$ to be non-trivial in $\mathcal{R}_k$, we need $k\geq 5$.
For any $k\geq 5$, we have $X_{{\bf m}'^{(k)}, {\bf m}} = Y_{{\bf m}'^{(k)}, {\bf m}} = \emptyset$, and $\mathrm{LC}_k(\m',\m)$ is satisfied. 

For $k=5$ or $6$, we have $X_{{\bf m}^{(k)}, {\bf m}'} = Y_{{\bf m}^{(k)}, {\bf m}'} = \emptyset$, so $\mathrm{LC}_k(\m',\m)$ is also satisfied and $[\overline{\Z}(\m)\times \overline{\Z}(\m')]= [\overline{\Z}(\m+\m')]$.

For $k\geq 7$, we have $X_{{\bf m}^{(k)}, {\bf m}'} = \{(1,1)\}$ and $Y_{{\bf m}^{(k)}, {\bf m}'} \subset \{(1,1)\}$, so $\mathrm{LC}_k(\m',\m)$ is not satisfied.

The $U_q(\widehat{\mathfrak{sl}_k})$-modules ($k \ge 5$) corresponding to ${\bf m} = [-4,-1]$ and ${\bf m}' = [-1,2]$ are $L(Y_{4,5})$ and $L(Y_{4,-1})$, respectively. Therefore for $k=5,6$, the $U_q(\widehat{\mathfrak{sl}_k})$ module $L(Y_{4,5}) \otimes L(Y_{4,-1})$ is simple. For $k \ge 7$, the $U_q(\widehat{\mathfrak{sl}_k})$ module $L(Y_{4,5}) \otimes L(Y_{4,-1})$ is not simple. In the case of $k=7$, we have short exact sequence
\begin{align*}
0 \to  L(Y_{1,2})  \to L(Y_{4,5}) \otimes L(Y_{4,-1}) \to L( Y_{4,5}Y_{4,-1} ) \to 0.
\end{align*}
In the case of $k \ge 8$, we have short exact sequence
\begin{align*}
0 \to L(Y_{1,2}) \otimes L(Y_{7,2})  \to L(Y_{4,5}) \otimes L(Y_{4,-1}) \to L( Y_{4,5}Y_{4,-1} ) \to 0.
\end{align*}

\end{example}

\begin{example}\label{example:2nd in main}
Let 
\begin{align*}
& {\bf m} = [-4, -3] + [-5, -4], \\
& {\bf n} = [0, 1] + [-1, 0]+ [-2, -2]  + [-2, -1] + [-3, -3] + [-3, -3] + [-5, -4]. 
\end{align*}
Then the corresponding monomials are 
\[
M_{\bf m} = Y_{2,7} Y_{2,9},  \quad M_{\bf n} = Y_{2,-1} Y_{2,1} Y_{1,4} Y_{2,3} Y_{1,6}^2 Y_{2,9}.
\] 
In $\mathcal{R}$, we have 
\begin{align*}
[\Z(\m)][\Z(\n)] = [\Z(\m+\n)] + [\Z({\m}')], 
\end{align*}
where ${\m}' = [0, 1] + [-1, 0] + [-4, -1] + [-2, -2] + [-5, -3] + [-3, -3] + [-5, -4]$.

For $k=3$, we have  
\begin{align*}
\chi_q(L(M_{\m})) \chi_q(L(M_{\n})) = \chi_q(L(M_{\m} M_{\n})).
\end{align*}
The tensor product $L(M_{\m}) \otimes L(M_{\n})$ is simple. Let us check the conditions ${\rm LC}_3(\m,\n)$ and ${\rm LC}_3(\n,\m)$. We have 

\begin{align*}
& X^{(3)}_{\m,\n} = \{(1, 3), (2, 5), (2, 6) \}, \quad 
 Y^{(3)}_{\m,\n} =  \{ (1, 5), (1, 6), (2, 7) \}.
\end{align*}
There is a matching function from $X^{(3)}_{\m,\n}$ to $Y^{(3)}_{\m,\n}$: $(1,3) \mapsto (1,5)$, $(2,5) \mapsto (2,7)$, $(2,6) \mapsto (1,6)$. Therefore ${\rm LC}_3(\m,\n)$ is satisfied. 

We have 
\begin{align*}
& X^{(3)}_{\n,\m} = \{(7,1) \}, \quad 
 Y^{(3)}_{\n,\m} =  \{ (7,2) \}.
\end{align*}
There is a matching function from $X^{(3)}_{\n,\m}$ to $Y^{(3)}_{\n,\m}$: $(7,1) \mapsto (7,2)$. Therefore ${\rm LC}_3(\n,\m)$ is satisfied. 

For $k=4$, we have  
\begin{align*}
\chi_q(L(M_{\m})) \chi_q(L(M_{\n})) = \chi_q(L(M_{\m} M_{\n})) + \chi_q(L(Y_{2,-1}Y_{2,1}Y_{1,4}Y_{3,8}Y_{1,6}Y_{2,9})).
\end{align*}
The tensor product $L(M_{\m}) \otimes L(M_{\n})$ is not simple. Let us check the conditions ${\rm LC}_4(\m,\n)$ and ${\rm LC}_4(\n,\m)$. We have 
\begin{align*}
& X^{(4)}_{\m,\n} = \{(1, 3), (1,4), (2, 5), (2, 6) \}, \quad 
 Y^{(4)}_{\m,\n} =  \{ (1, 5), (1, 6), (2, 7) \}.
\end{align*}
There is no matching function from $X^{(4)}_{\m,\n}$ to $Y^{(4)}_{\m,\n}$. Therefore ${\rm LC}_4(\m,\n)$ is not satisfied. 
 
We have 
\begin{align*}
& X^{(4)}_{\n,\m} = \{(7,1) \}, \quad 
 Y^{(4)}_{\n,\m} =  \{ (7,1), (7,2) \}.
\end{align*}
There is a matching function from $X^{(4)}_{\n,\m}$ to $Y^{(4)}_{\n,\m}$: $(7,1) \mapsto (7,2)$. Therefore ${\rm LC}_4(\n,\m)$ is satisfied.

For $k \ge 5$, we have that 
\begin{align*}
\chi_q(L(M_{\m})) \chi_q(L(M_{\n})) = \chi_q(L(M_{\m} M_{\n})) + \chi_q(L(Y_{2,-1}Y_{2,1}Y_{4,5}Y_{1,4}Y_{3,8}Y_{1,6}Y_{2,9})).
\end{align*}
The tensor product $L(M_{\m}) \otimes L(M_{\n})$ is not simple.
We have 
\begin{align*}
& X^{(k)}_{\m,\n} = \{(1, 3), (1,4), (2, 5), (2, 6) \}, \quad 
 Y^{(k)}_{\m,\n} =  \{ (1, 5), (1, 6), (2, 7) \}.
\end{align*}
There is no matching function from $X^{(k)}_{\m,\n}$ to $Y^{(k)}_{\m,\n}$. Therefore ${\rm LC}_k(\m,\n)$ is not satisfied. 

We have 
\begin{align*}
& X^{(k)}_{\n,\m} = \{(7,1) \}, \quad 
 Y^{(k)}_{\n,\m} =  \{ (7,1), (7,2) \}.
\end{align*}
There is a matching function from $X^{(k)}_{\n,\m}$ to $Y^{(k)}_{\n,\m}$: $(7,1) \mapsto (7,2)$. Therefore ${\rm LC}_k(\n,\m)$ is satisfied.  
\end{example}

\section{Proofs of main results}\label{sect_proof} 

In this section, we mostly drop the subscripts and write $\m$, $\n$, as we only consider multisegments.

\subsection{Proof of Proposition \ref{proposition_LC_LI_cusp}}\label{subsec:the case that m is a segment with length 1}




Suppose here that $\m$ is of length 1 and write $\m=[a,a]$. Recall that $k\in\mathbb{Z}_{>0}$ is fixed.

\begin{remark} Note that in this case, the corresponding tableaux and monomial are $T=T_{a, a}$ ($T_{a,b}$ is defined in Section \ref{def:weakly separted k-subsets}) and $M_T = Y_{1,-2a}$. 
\end{remark}

Let $\n$ be any multisegment classifying an irreducible object in $\mathcal{C}_k$, then $\n=\Delta_1 + \cdots + \Delta_N$, where all $\Delta_i$ have length $e(\Delta_i)-b(\Delta_i) +1<k$.

As $X_{\m,\n}\subset \{1\}\times [N]$, by abuse of notations, we consider $X_{\m,\n}$ as a subset of $[N]$. Then 
\[ X_{\m,\n} = \left\lbrace i \in [N] \mid b(\Delta_i) = a+1 \right\rbrace, \quad X_{\m,\n}^{(k)} = \left\lbrace i \in [N] \mid b(\Delta_i) = a+1, e(\Delta_i)<k+a \right\rbrace.\]

It is clear that $X_{\m,\n}^{(k)} \subset X_{\m,\n}$. Suppose that $i \in X_{\m,\n}^{(k)}$. Then $b(\Delta_i) = a+1$ for each $i \in [N]$. By the condition on the length of the segments in $\n$, we have for each $i \in [N]$, $e(\Delta_i) - b(\Delta_i)+1<k$. Therefore $e(\Delta_i) - (a+1) + 1<k$ and hence $e(\Delta_i)<k+a$. Thus $i \in X_{\m,\n}^{(k)}$. It follows that
\[  X_{\m,\n} = X_{\m,\n}^{(k)}.\]
Similarly, 
\[ Y_{\m,\n} = Y_{\m,\n}^{(k)}.\]

Thus in this case, the condition ${\rm LC}(\m,\n)$ is equivalent to the condition ${\rm LC}_k(\m,\n)$.

We need to prove the following result, which combined with Proposition~\ref{proposition_LI_implies_LIk}, proves Proposition~\ref{proposition_LC_LI_cusp} in this case.

\begin{proposition}
For $\rho =\Z([a,a])$ and $\pi = \Z(\n)$, where $\n \in \Irr_k$, 
\[{\rm LI}_k(\overline{\rho},\overline{\pi}) \Rightarrow {\rm LI}(\rho,\pi).\]
\end{proposition}

We use the following result, stated in \cite[Theorem 5.11]{LM16}, based on results appearing in \cite[Theorem 7.5]{Min09} and \cite[Theorem 2.2.1]{Jan07}.

\begin{proposition}\label{proposition_soc_rho_Z(n)}
Let $f$ be the best matching between $X_{[a,a],\n}$ and $Y_{[a,a],\n}$. Then
\begin{enumerate}
	\item ${\rm soc}(\rho\times {\rm Z}(\n)) = {\rm Z}(\n+ [a,a])$ if and only if $f$ is a function from $X_{[a,a],\n}$ to $Y_{[a,a],\n}$.
	\item If $f$ is a not function from $X_{[a,a],\n}$ to $Y_{[a,a],\n}$, and $i\in X_{[a,a],\n}$ is the minimal index which does not belong to the domain of $f$, then
	\[{\rm soc}(\rho\times {\rm Z}(\n)) = {\rm Z}(\n - \Delta_i + [a,e(\Delta_i)]).\]
\end{enumerate}
\end{proposition}

\begin{proof}
From Proposition~\ref{proposition_soc_rho_Z(n)}, if the condition ${\rm LI}(\rho,\pi)$ is not satisfied, then ${\rm soc}(\rho\times {\rm Z}(\n)) = {\rm Z}(\n - \Delta_i + [a,e(\Delta_i)])$. However, for $i\in X_{[a,a],\n}$, $e(\Delta_i)<k+a$, so the segment $[a,e(\Delta_i)]$ has length strictly lower than $k+1$. Hence, the image under $\overline{\mathcal{Q}}$ of ${\rm Z}(\n - \Delta_i + [a,e(\Delta_i)])$ is a simple object in $\mathcal{C}_k$, which is different from $\overline{\rm Z}(\n)$. Thus, the condition ${\rm LI}_k(\overline{\rho},\overline{\pi})$ is also not satisfied.
\end{proof}

\begin{example}
(1) Let $k \ge 3$, $\rho = \oZ([0,0])$, ${\bf n} = \Delta_1+\Delta_2$, $\Delta_1=[1,2]$, $\Delta_2=[0,1]$. Then $X^{(k)}_{\rho, {\bf n}} = \{1\}$, $Y^{(k)}_{\rho, {\bf n}}=\{2\}$. The best $\rightsquigarrow_k$-matching $f$ between $X^{(k)}_{\rho, {\bf n}}$ and $Y^{(k)}_{\rho, {\bf n}}$ is given by $f(1)=2$. The domain of $f$ is $X^{(k)}_{\rho, {\bf n}}$ and hence $f$ is a function from $X^{(k)}_{\rho, {\bf n}}$ to $Y^{(k)}_{\rho, {\bf n}}$. We have that ${\rm soc}(\rho \times \oZ({\bf n})) = \oZ([0,0]+{\bf n})$. In the language of quantum affine algebra, this means for $U_q(\widehat{\mathfrak{sl}_k})$-modules $L(Y_{1,0})$ and $L(Y_{2,-3}Y_{2,-1})$, the socle of $L(Y_{1,0}) \otimes L(Y_{2,-3}Y_{2,-1})$ is $L(Y_{1,0}Y_{2,-3}Y_{2,-1})$.

(2) Let $k=4$, $\rho = \oZ([0,0])$, ${\bf n} = \Delta_1+\Delta_2$, $\Delta_1=[2,3]$, $\Delta_2=[1,2]$. Then $X^{(k)}_{\rho, {\bf n}} = \{2\}$, $Y^{(k)}_{\rho, {\bf n}}=\emptyset$. The domain of the best $\rightsquigarrow_k$-matching $f$ between $X^{(k)}_{\rho, {\bf n}}$ and $Y^{(k)}_{\rho, {\bf n}}$ is $\emptyset$. Therefore $f$ is not a function from $X^{(k)}_{\rho, {\bf n}}$ to $Y^{(k)}_{\rho, {\bf n}}$. We have that ${\rm soc}(\rho \times \oZ({\bf n})) = \oZ({\bf n} - \Delta_2 + {}^{-}\! \Delta_2)=\oZ([2,3]+[0,2])$ in $\mathcal{C}_k$. In the language of quantum affine algebra, this means for $U_q(\widehat{\mathfrak{sl}_4})$-modules $L(Y_{1,0})$ and $L(Y_{2,-5}Y_{2,-3})$, the socle of $L(Y_{1,0}) \otimes L(Y_{2,-5}Y_{2,-3})$ is $L(Y_{2,-5}Y_{3,-2})$.

Let $k=3$, $\rho = \oZ([0,0])$, ${\bf n} = \Delta_1+\Delta_2$, $\Delta_1=[2,3]$, $\Delta_2=[1,2]$. Then $X^{(k)}_{\rho, {\bf n}} = \{2\}$, $Y^{(k)}_{\rho, {\bf n}}=\emptyset$. The domain of the best $\rightsquigarrow_k$-matching $f$ between $X^{(k)}_{\rho, {\bf n}}$ and $Y^{(k)}_{\rho, {\bf n}}$ is $\emptyset$. Therefore $f$ is not a function from $X^{(k)}_{\rho, {\bf n}}$ to $Y^{(k)}_{\rho, {\bf n}}$. We have that ${\rm soc}(\rho \times \oZ({\bf n})) = \oZ({\bf n} - \Delta_2 + {}^{-}\! \Delta_2)=\oZ([2,3]+[0,2])=\oZ([2,3])$ in $\mathcal{C}_k$. In the language of quantum affine algebra, this means for $U_q(\widehat{\mathfrak{sl}_3})$-modules $L(Y_{1,0})$ and $L(Y_{2,-5}Y_{2,-3})$, the socle of $L(Y_{1,0}) \otimes L(Y_{2,-5}Y_{2,-3})$ is $L(Y_{2,-5})$.
\end{example}

\subsection{Proof of Proposition~\ref{prop: ladder times ladder case}} \label{sec:proof of the case of ladder times ladder}

First we recall results in \cite{Gur19}, \cite{Gur21}. 
%
For $\lambda=(\lambda_1, \ldots, \lambda_m), \mu=(\mu_1, \ldots, \mu_m) \in \ZZ^m$, denote
\[
\m(\lambda, \mu) = \sum_{i=1}^m [\lambda_i, \mu_{i}].
\] 
For $x$ a permutation in $S_m$ such that $\lambda_i \leq \mu_{x(i)} +1$, for all $1\leq i\leq m$, 
\begin{equation}
\m_x = \m_x(\lambda, \mu) = \sum_{i=1}^m [\lambda_i, \mu_{x(i)}].
\end{equation}

For $\lambda=(\lambda_1 \le \ldots \le \lambda_m), \mu=(\mu_1 \le \ldots \le \mu_m) \in \ZZ^m$, let $S_m^{(\lambda,\mu)}$ be the double-coset of permutations inside $S_m$ up to the following equivalence
\begin{equation}\label{eq:equivalence relation for permutations}
x \sim x' \quad \Leftrightarrow \quad \m_x(\lambda,\mu) = \m_{x'}(\lambda,\mu).
\end{equation}

Let $\m, \n$ be multisegments and let $\lambda=(\lambda_1 \le \ldots \le \lambda_m), \mu=(\mu_1 \le \ldots \le \mu_m) \in \ZZ^m$ be the tuples of integers which are the beginning and ending of the multisegments in $\m+\n$, respectively.

For $1\leq i,j \leq m$ such that $\lambda_i \leq \mu_j+1$, consider the segment $\Delta = [\lambda_i,\mu_j]$, and define a sequence 
\begin{align*}
{\rm Seq}(\m,\n,\Delta) = ((\Delta_1', \Delta_{i_1}, \n_1'), (\Delta_2', \Delta_{i_2}, \n_2'), \ldots, (\Delta_r', \Delta_{i_r}, \n_r'))
\end{align*}
as follows: 
\begin{itemize}
\item First we order the segments in $\m + \n = \sum_{i} \Delta_i$ such that $b(\Delta_i )\leq b(\Delta_{i+1})$, and if $b(\Delta_i)=b(\Delta_{i+1})$, then $e(\Delta_i)\geq e(\Delta_{i+1})$.

\item Step 1. Take the smallest $i_1$ such that $b(\Delta_{i_1}) = b(\Delta)$. Let $\Delta_1'$ be the sub-segment of $\Delta_{i_1}$ whose support is the intersection of the supports of $\Delta_{i_1}$ and $\Delta$. If $\Delta_{i_1}$ is a segment of $\m$, then $\n_1'=\m$. Otherwise, $\n_1'=\n$.

\item Step 2. Now consider the smallest $i_2$ such that the support of $\Delta_{i_2}$ intersects the support of $\Delta \setminus \Delta_1'$, and $\Delta_{i_2}$ and $\Delta_{i_1}$ (the $\Delta_{i_1}$ in previous step) are neither both in $\m$ nor both in $\n$. Let $\Delta_2'$ be the sub-segment of $\Delta_{i_2}$ whose support is the intersection of the supports of $\Delta_{i_2}$ and $\Delta \setminus \Delta_1'$. If $\Delta_{i_2}$ is a segment of $\m$, then $\n_2'=\m$. Otherwise, $\n_2'=\n$.

\item Continue this procedure, at Step $j$, take the smallest $i_j$ such that the support of $\Delta_{i_j}$ intersects the support of $\Delta \setminus ( \cup_{s=1}^{j-1} \Delta_s' )$, and $\Delta_{i_j}$ and $\Delta_{i_{j-1}}$ (the $\Delta_{i_{j-1}}$ in previous step) are neither both in $\m$ nor both in $\n$. Let $\Delta_j'$ be the sub-segment of $\Delta_{i_j}$ whose support is the intersection of supports of $\Delta_{i_j}$ and $\Delta \setminus ( \cup_{s=1}^{j-1} \Delta_s' )$. If $\Delta_{i_j}$ is a segment of $\m$, then $\n_j'=\m$. Otherwise, $\n_j'=\n$.

\item Continue this procedure until at some step $j$, there is no $\Delta_{i_j}$ whose support intersects the support of $\Delta \setminus ( \cup_{s=1}^{j-1} \Delta_s' )$, and $\Delta_{i_j}$ and $\Delta_{i_{j-1}}$ (the $\Delta_{i_{j-1}}$ in previous step) are neither both in $\m$ nor both in $\n$.
\end{itemize}

The multisegments $\m, \n$ are said to be \emph{tiled} by $\Delta$ if the union of the supports of $\Delta_j'$, $j=1,\ldots,r$, is the support of $\Delta$, and for all $j$, 
$e(\Delta_{j}')=e(\Delta_{i_j})$. 

We say that removing ${\rm Seq}(\m,\n,\Delta)$ from $\m, \n$ is to remove $\Delta_j'$ from $\Delta_{i_j}$ for each $j=1,\ldots, r$. After removing, we obtain two multisegments $\m', \n'$ whose segments are sub-segments of $\m, \n$ respectively. Moreover, there exists $x\in S_m$ such that $\Delta + \m' +\n' =\m_{x}(\lambda,\mu)$.

\begin{example}
Let $\m=[-6,-1] + [-2,3] + [-1,4]$, $\n = [-4,1] + [0,2]$. Then $\m+\n = \sum_{i} \Delta_i = [-6,-1] + [-4,1] + [-2,3] + [-1,4] + [0,2]$ (ordered using $\le_b$).  For $\Delta = [-4, 3]$, 
\[
{\rm Seq}(\m, \n, \Delta) =( ( ([-4, 1], \Delta_2, \n), ([2,3], \Delta_3, \m) ) ).
\] 
In this case, the multisegments $\m, \n$ are tiled by $\Delta$.

For $\Delta=[-6,2]$, 
\[
{\rm Seq}(\m, \n, \Delta) = (( [-6, -1], \Delta_1, \m), ([0,1], \Delta_2, \n), ([2,2], \Delta_3, \m) ).
\]
In this case, the multisegments $\m, \n$ are not tiled by $\Delta$ since the end point of $[2,2]$ and $\Delta_3=[-2,3]$ are different. 
\end{example}


Now, write the segments of $\m_x = \sum_{i=1}^r \Delta_i$ such that $b(\Delta_i )\leq b(\Delta_{i+1})$, and if $b(\Delta_i)=b(\Delta_{i+1})$, then $e(\Delta_i)\geq e(\Delta_{i+1})$. Using the terminology of Gurevich \cite{Gur21}, we say that $\m_x$ \emph{tiles} $(\m, \n)$ if the following procedure stops only after step $r$. 
\begin{itemize}
\item Step 1. If the multisegments $\m, \n$ are tiled by $\Delta_1$, then let $\m_1, \n_1$ be the multisegments obtained by removing ${\rm Seq}(\m,\n,\Delta_1)$ from $\m, \n$. 

\item Step 2. If the multisegments $\m_1, \n_1$ are tiled by $\Delta_2$, then let $\m_2, \n_2$ be the multisegments obtained by removing ${\rm Seq}(\m_1,\n_1,\Delta_2)$ from $\m_1, \n_1$. 

\item Step $j$. Continue this procedure, if the multisegments $\m_{j-1}, \n_{j-1}$ are tiled by $\Delta_j$, then let $\m_j, \n_j$ be the multisegments obtained by removing ${\rm Seq}(\m_{j-1},\n_{j-1},\Delta_j)$ from $\m_{j-1}, \n_{j-1}$. 
\end{itemize}



A permutation $x\in S_m$ is \emph{321-avoiding} if there is no sequence $1\leq i_1<i_2<i_3\leq m$ such that $x(i_3)<x(i_2)<x(i_2)$. It is well-known that a permutation $x\in S_m$ is 321-avoiding if and only if there exists a disjoint partition $I\sqcup J = \{1,2,\ldots,m\}$ such that $x_{|I}$ and $x_{|J}$ are strictly increasing. For completeness, we prove the following similar result.

\begin{lemma}\label{lem: 321 not sum of ladders}
    The multisegment $\m_{x}(\lambda,\mu)$ is a sum of two ladders if and only if the longest representative $x'$ of $x$ in $S_m^{(\lambda,\mu)}$ is 321-avoiding.
\end{lemma}

\begin{proof}
If $[\lambda_i,\mu_{x(i)}]$ and $[\lambda_j,\mu_{x(j)}]$ form a ladder, then (without loss of generality) $\lambda_i < \lambda_j$ and $\mu_{x(i)} < \mu_{x(j)}$. So $i<j$ and $x(i)<x(j)$. This proves the direct implication.

Now suppose that $\m_{x}(\lambda,\mu)$ is not the sum of two ladders. Then there exists $i<j<l$ such that $\lambda_i \leq \lambda_j \leq \lambda_l$ and $\mu_{x(l)} \leq \mu_{x(j)} \leq \mu_{x(i)}$. Thus if $x'$ is the longest representative of $x$ in $S_m^{(\lambda,\mu)}$, we have $x'(l)<x'(j)<x'(i)$.
\end{proof}

We use the following result, which is the main result of \cite{Gur21}.

\begin{theorem}[{\cite[Theorem 1.2]{Gur19}, \cite[Theorem 1.1]{Gur21}}]
Let $\m, \n$ be ladders in ${\rm Mult}_k$ and let $\lambda=(\lambda_1 \le \ldots \le \lambda_m), \mu=(\mu_1 \le \ldots \le \mu_m) \in \ZZ^m$ be the tuples of integers which are the beginning and ending of $\m+\n$ respectively. Then we have
\begin{align} \label{eq:Zm times Zn ladders in R}
[\Z(\m) \times \Z(\n)] = \sum_{x \in S(\m, \n)} [\Z(\m_x)] \in \mathcal{R},
\end{align}
where $S(\m, \n) \subset S_m^{(\lambda,\mu)}$ is formed of permutations $x \in S_m$ such that $\lambda_i \leq \mu_{x(i)} +1, \forall i $, $\m_{x}$ tiles $\m, \n$ and the longest representative $x'$ of $x$ in $S_m^{(\lambda,\mu)}$ is 321-avoiding.
\end{theorem}

Therefore in $\mathcal{R}_k$, for ladders  $\m, \n$ in ${\rm Mult}_k$,  we have 
\begin{align} \label{eq:Zm times Zn ladders in Rk}
[\oZ(\m) \times \oZ(\n)] = \sum_{x \in S(\m, \n)} [\oZ(\om_x)] \in \mathcal{R}_k.
\end{align}

\begin{remark}
In \eqref{eq:Zm times Zn ladders in R}, all terms corresponding to permutations $x\in S_m$ for which there is $i$ such that $\lambda_i > \mu_{x(i)} +1$ are sent to 0. In \eqref{eq:Zm times Zn ladders in Rk}, additionally, all terms corresponding to permutations $x\in S_m$ for which there is $i$ such that $\mu_{x(i)} >\lambda_i + k - 1$ are sent to 0.
\end{remark}

For $\m,\n\in{\rm Mult}_k$ ladders, write $\m = \Delta_1 + \Delta_2 + \cdots + \Delta_N$, and $\n = \Delta'_1 + \Delta'_2 + \cdots + \Delta'_{N'}$, in their \emph{left} aligned form ($b(\Delta_i) <b(\Delta_{i+1})$ and $b(\Delta'_i) <b(\Delta'_{i+1})$). Let $X= X_{\m,\n}^{(k)}$ and $Y=Y_{\m,\n}^{(k)}$, and denote by ${\rm NC}_k(\m,\n)$ the condition that there exists indices $1\leq i\leq N$, $1\leq j\leq N'$, and $r\geq 0$ such that
\begin{itemize}
	\item $(i+l, j+l)\in X$, for all $0\leq l \leq r$,
	\item $(i, j-1)$ and $(i+r+1, j+r)$ are not in $Y$,
	\item $(i+l+1,j+l) \in Y\setminus X$, for all $0\leq l \leq r-1$.
\end{itemize}
Note that $(i,j) \in Y\setminus X$ implies that $b(\Delta_i) = b(\Delta'_j)$ or $e(\Delta_i) = e(\Delta'_j)$ .

The following is easy an modification from \cite[Lemma 6.21]{LM16}:
\begin{lemma}\label{lemma_LC_NC}
For two ladders $\m,\n \in {\rm Mult}_k$, the negation of ${\rm LC}_k(\m,\n)$ is equivalent to ${\rm NC}_k(\m,\n)$.
\end{lemma}

We prove Proposition~\ref{prop: ladder times ladder case}, via the following equivalent result.

\begin{proposition} \label{prop:ladder times ladder}
For ladders $\m, \n$ in ${\rm Mult}_k$, the conditions ${\rm LC}_k(\m, \n)$ and ${\rm LC}_k(\n, \m)$ hold if and only if there is no $\m_x$ appearing on the right hand side of (\ref{eq:Zm times Zn ladders in Rk}) such that $\m_x \ne \m+\n$ and $\om_x$ is not $0$.
\end{proposition}

\begin{proof}
Throughout this proof, multisegments $\sum_i \Delta_i$ will be \emph{ordered} as follows: $b(\Delta_i)\leq b(\Delta_{i+1})$, and if $b(\Delta_i)=b(\Delta_{i+1})$, then $e(\Delta_i)\geq e(\Delta_{i+1})$ (note that this order is different from $\leq_b$, but more adapted to look for the longest representative associated to the permutation).

We write $\m=\sum_{i} \Delta_i^{(\m)}$, $\n=\sum_{i} \Delta_i^{(\n)}$ using the order above. Let $\lambda=(\lambda_1 \le \ldots \le \lambda_m), \mu=(\mu_1 \le \ldots \le \mu_m) \in \ZZ^m$ be the tuples of integers which are the beginning and ending of $\m+\n$ respectively. We write $\m+\n=\sum_j \Delta_j$ using the same order. We have $\m+\n=\m_{x}(\lambda, \mu)$ for some $x\in S_m$.

$( \Leftarrow )$ Suppose, without loss of generality, that the condition ${\rm LC}_k(\m, \n)$ does not hold. From Lemma~\ref{lemma_LC_NC}, we know that the condition ${\rm NC}_k(\m, \n)$ holds. Let $1\leq i\leq N$, $1\leq j\leq N'$, and $r\geq 0$ be the data associated to the condition ${\rm NC}_k(\m, \n)$.

Let $\m'$ be the multi-segment obtained from $\m + \n$ by replacing, for all $0\leq l\leq r$, the segments $\Delta^{(\m)}_{i+l}, \Delta^{(\n)}_{j+l}$ by $\left(\Delta^{(\m)}_{i+l}\cup\Delta^{(\n)}_{j+l}\right) , \left(\Delta^{(\m)}_{i+l}\cap\Delta^{(\n)}_{j+l}\right)$. Since $\Delta^{(\m)}_{i+l}\prec_k\Delta^{(\n)}_{j+l}$, the resulting multisegment has the same set of extremities as $\m+\n$, thus there exists $x'\in S_m$ such that $\m' = \m_{x'}(\lambda,\mu)$. Moreover, $[\overline{\Z}(\m_{x'})] \neq 0 \in \mathcal{R}_k$.


Moreover, for all $0\leq l \leq r$, the tiling sequences of $\left(\Delta^{(\m)}_{i+l}\cup\Delta^{(\n)}_{j+l}\right)$ and $\left(\Delta^{(\m)}_{i+l}\cap\Delta^{(\n)}_{j+l}\right)$ are:
\begin{align*}
    {\rm Seq}\left(\Delta^{(\m)}_{i+l}\cup\Delta^{(\n)}_{j+l}\right) & = \left((\Delta^{(\m)}_{i+l}, \m),\left(\Delta^{(\n)}_{j+l}\setminus \left( \Delta^{(\n)}_{j+l} \cap \Delta^{(\m)}_{i+l}\right),\n\right)\right),\\
    {\rm Seq}\left(\Delta^{(\m)}_{i+l}\cap\Delta^{(\n)}_{j+l}\right) & = \left(\Delta^{(\m)}_{i+l}\cap\Delta^{(\n)}_{j+l},\n\right).
\end{align*}
Indeed, since $(i+l,j+l-1)\in Y\setminus X$ (or if $l=0$, $(i,j-1)\notin Y$), so there is no $\Delta^{(\n)}_{j'}\in \n$ such that $e(\Delta^{(\m)}_{i+l})< e(\Delta^{(\n)}_{j'}) < e(\Delta^{(\n)}_{j+l})$. Hence $\m_{x'}$ tiles $\m,\n$.

Finally, in order to prove that the longest representative of $x'$ in $S_m^{(\lambda, \mu)}$ is 321-avoiding, thanks to Lemma~\ref{lem: 321 not sum of ladders}, we only have to prove that $\m'$ is the sum of two ladders. 
We have 
\[
\m' = \sum_{l<i, l>i+r}\Delta^{(\m)}_{l} + \sum_{l<j, l>j+r}\Delta^{(\n)}_{l} + \sum_{l=0}^r\left( \Delta^{(\m)}_{i+l}\cup\Delta^{(\n)}_{j+l} + \Delta^{(\m)}_{i+l}\cap\Delta^{(\n)}_{j+l}\right).
\]

Suppose there exists $\Delta^{(\n)}_{j-1} \in \n$, then by ${\rm NC}_k(\m, \n)$, we have $(i,j-1) \notin Y$. Then, either
\begin{center}
    \begin{tabular}{ll}
  $(a)$ $b(\Delta^{(\n)}_{j-1})< b(\Delta^{(\m)}_{i})$   & or  $(b)$ $e(\Delta^{(\m)}_{i})< b(\Delta^{(\n)}_{j-1})$,\\
  or $(c)$ $e(\Delta^{(\n)}_{j-1})< e(\Delta^{(\m)}_{i})$   & or  $(d)$ $b(\Delta^{(\m)}_{i}) + k -2 < e(\Delta^{(\n)}_{j-1})$.
\end{tabular}
\end{center}
Cases $(b)$ and $(d)$ are impossible, since $\Delta^{(\m)}_{i}\prec_k \Delta^{(\n)}_{j}$, and $\Delta^{(\n)}_{j-1}$ and $\Delta^{(\n)}_{j}$ form a ladder.

\begin{center}
\begin{tikzpicture}[scale=0.7]
\draw[ultra thick,blue] (0,1.5) -- (5,1.5);
\draw[ultra thick,blue, dotted] (5,1.5) -- (7,1.5);
\draw[ultra thick] (2,1) -- (6,1);
\draw[ultra thick,red] (3,0.5) -- (5,0.5);
\draw[ultra thick, red, dotted ] (0,0.5) -- (3,0.5);
\node at (1.3,1) {\small $\Delta^{(\m)}_{i}$};
\node[blue] at (-0.8,1.5) {\small $\Delta^{(\m)}_{j-1}$};
\node[red] at (-0.8,0.5) {\small $\Delta^{(\m)}_{j-1}$};
\node[blue] at (7.8,1.5) {$(a)$};
\node[red] at (5.8,0.5) {$(c)$};
\end{tikzpicture}
\end{center}

In case $(a)$, since we also have $e(\Delta^{(\n)}_{j-1})< e(\Delta^{(\n)}_{j}) = e(\Delta^{(\m)}_{i}\cup\Delta^{(\n)}_{j})$, then 
\[
  \m_1 =\sum_{l<j}\Delta^{(\n)}_{l} + \sum_{l=0}^r \Delta^{(\m)}_{i+l}\cup\Delta^{(\n)}_{j+l}
\]
forms a ladder.
Additionally, if there exists $\Delta^{(\m)}_{i-1} \in \m$, then we also have $b(\Delta^{(\m)}_{i-1})< b(\Delta^{(\m)}_{i}) < b(\Delta^{(\n)}_{j}) =b(\Delta^{(\m)}_{i}\cap\Delta^{(\n)}_{j})$, and $e(\Delta^{(\m)}_{i-1}) < e(\Delta^{(\m)}_{i})=e(\Delta^{(\m)}_{i}\cap\Delta^{(\n)}_{j})$. 
Then 
\[
 \m_2 = \sum_{l<i}\Delta^{(\m)}_{l} + \sum_{l=0}^r \Delta^{(\m)}_{i+l}\cap\Delta^{(\n)}_{j+l}
\]
forms a ladder.

In case $(c)$, since we also have $b(\Delta^{(\n)}_{j-1})< b(\Delta^{(\n)}_{j}) = b(\Delta^{(\m)}_{i}\cap\Delta^{(\n)}_{j})$, then 
\[
 \m_2 =\sum_{l<j}\Delta^{(\n)}_{l} + \sum_{l=0}^r \Delta^{(\m)}_{i+l}\cap\Delta^{(\n)}_{j+l}
\]
forms a ladder.
Additionally, if there exists $\Delta^{(\m)}_{i-1} \in \m$, then we also have $e(\Delta^{(\m)}_{i-1})< e(\Delta^{(\m)}_{i}) < e(\Delta^{(\n)}_{j}) =e(\Delta^{(\m)}_{i}\cup\Delta^{(\n)}_{j})$, and $b(\Delta^{(\m)}_{i-1}) < b(\Delta^{(\m)}_{i})=b(\Delta^{(\m)}_{i}\cup\Delta^{(\n)}_{j})$. 
Then 
\[
  \m_1 =\sum_{l<i}\Delta^{(\m)}_{l} + \sum_{l=0}^r \Delta^{(\m)}_{i+l}\cup\Delta^{(\n)}_{j+l}
\]
forms a ladder.

If there exists $\Delta^{(\m)}_{i+r=1} \in \m$, then by ${\rm NC}_k(\m, \n)$, we have $(i+r+1,j+r) \notin Y$. Then by a distinction of cases as above, we can add the segments in $\sum_{l>i+r}\Delta^{(\m)}_{l} $ and $\sum_{l>j+r}\Delta^{(\n)}_{l} $ to either $\m_1$ or $\m_2$ to order to build two ladders $\m_1',\m_2'$ such that $\m'=\m_1' + \m_2'$.

$( \Rightarrow )$ Suppose we have a 321-avoiding permutation $x'\in S_m^{(\lambda,\mu)}$ such that $\overline{\m}_{x'}\neq  \m + \n, 0$ and $\m_{x'}$ tiles $\m + \n$. 

Since $\overline{\m}_{x'}\neq  \m + \n$, there is at least one segment in $\overline{\m}_{x'}$ whose sequence ${\rm Seq}$ in the algorithm above is of length $\geq 2$. Take $\Delta$ to be the minimal such segment (for the order above). Replace $\m$ and $\n$ by the multisegments obtained after removing all the segments before $\Delta$ in $\overline{\m}_{x'}$, as they are all equal to segments of $\m + \n$, and denote the resulting multisegments by $\m$ and $\n$ again (in other words, we can assume that $\Delta$ is the first segment of $\overline{\m}_{x'}$). Now, let us assume without loss of generality that $\mathrm{Seq}(\m,\n,\Delta)$ starts as follows 
\[ \mathrm{Seq}(\m,\n,\Delta) = ((\Delta_1',\Delta_i^{(\m)}, \m) , (\Delta_2',\Delta_j^{(\n)}, \n), \ldots ).\]
As above, we use the notations $X=X_{\m,\n}^{(k)}$ and $Y=Y_{\m,\n}^{(k)}$.
As $\m_{x'}$ tiles $\m+\n$, necessarily $\Delta_1'= \Delta_i^{(\m)}$ and $b(\Delta) = b(\Delta_i^{(\m)})$. Moreover, $b(\Delta_j^{(\n)}) > b(\Delta_i^{(\m)})$, by definition of $\Delta_i^{(\m)}$, and $[e(\Delta_i^{(\m)}) +1 , e(\Delta)] \cap \Delta_j^{(\n)} \neq \emptyset$, thus $e(\Delta_j^{(\n)}) \geq e(\Delta_i^{(\m)}) +1$. Next, we distinguish 2 cases (as illustrated below).
Either, (A) : $b(\Delta_2')= e(\Delta_i^{(\m)}) +1$, or (B) : $b(\Delta_2')> e(\Delta_i^{(\m)}) +1$.

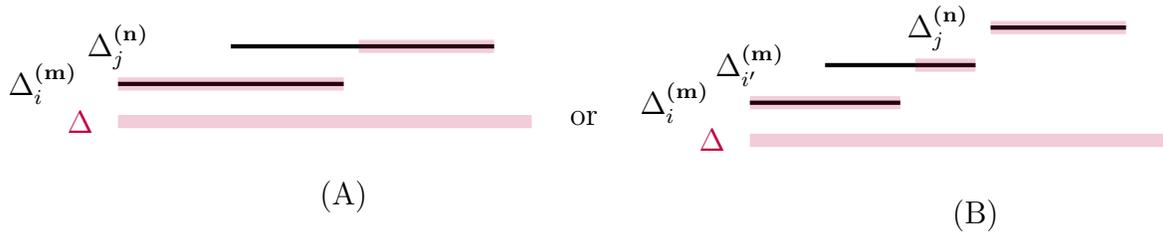
\begin{center}
\begin{tikzpicture}
\node at (0,0) {\begin{tikzpicture}
\draw[line width=5pt, purple, opacity=0.2] (0,0) -- (5.5,0);
\draw[ultra thick] (0,0.5) -- (3,0.5);
\draw[ultra thick] (1.5,1) -- (5,1);
\draw[line width=5pt, purple, opacity=0.2] (0,0.5) -- (3,0.5);
\draw[line width=5pt, purple, opacity=0.2] (3.2,1) -- (5,1);
\node[purple] at (-0.5,0) {$\Delta$};
\node at (-1,0.5) {$\Delta_i^{(\m)}$};
\node at (0,1) {$\Delta_j^{(\n)}$};
\node at (3,-1) {(A)};
\end{tikzpicture}};
\node at (4.2,0) {or};
\node at (8.4,0) {\begin{tikzpicture}
\draw[line width=5pt, purple, opacity=0.2] (0,0) -- (5.5,0);
\draw[ultra thick] (0,0.5) -- (2,0.5);
\draw[ultra thick] (1,1) -- (3,1);
\draw[ultra thick] (3.2,1.5) -- (5,1.5);
\draw[line width=5pt, purple, opacity=0.2] (0,0.5) -- (2,0.5);
\draw[line width=5pt, purple, opacity=0.2] (2.2,1) -- (3,1);
\draw[line width=5pt, purple, opacity=0.2] (3.2,1.5) -- (5,1.5);
\node[purple] at (-0.5,0) {$\Delta$};
\node at (-1,0.5) {$\Delta_i^{(\m)}$};
\node at (0,1) {$\Delta_{i'}^{(\m)}$};
\node at (2.5,1.5) {$\Delta_j^{(\n)}$};
\node at (3,-1) {(B)};
\end{tikzpicture}};
\end{tikzpicture}
\end{center}

In the case (A), $\Delta_i^{(\m)}\prec_k \Delta_j^{(\n)}$, thus $(i,j)\in X$.
Suppose there exists $\Delta_{j-1}^{(\n)} \in \n$, and $(i,j-1)\in Y$. Then since $\Delta$ is the first segment of $\overline{\m}_{x'}$, and by definition of $\Delta_{j}^{(\n)}$, $e(\Delta_i^{(\m)}) = e(\Delta_{j-1}^{(\n)})$ and $b(\Delta_i^{(\m)}) \leq b(\Delta_{j-1}^{(\n)})$. 

\begin{center}
    \begin{tikzpicture}
\draw[ultra thick] (0,0.5) -- (3,0.5);
\draw[ultra thick] (1.5,1.5) -- (5,1.5);
\draw[ultra thick] (0.5,1) -- (3,1);
\draw[line width=5pt, purple, opacity=0.2] (0,0.5) -- (3,0.5);
\draw[line width=5pt, purple, opacity=0.2] (3.2,1.5) -- (5,1.5);
\node at (-0.7,0.5) {$\Delta_i^{(\m)}$};
\node at (0,1) {$\Delta_{j-1}^{(\n)}$};
\node at (0.8,1.5) {$\Delta_j^{(\n)}$};
\end{tikzpicture}
\end{center}

Thus there exists $\Delta', \Delta'' \in \overline{\m}_{x'}$ such that $b(\Delta'),b(\Delta'') \geq b(\Delta) = b(\Delta_i^{(\m)})$ and $e(\Delta')=e(\Delta'') = e(\Delta_i^{(\m)}) < e(\Delta)$. Necessarily, the segments $\Delta,\Delta'$ and $\Delta''$ are pairwise in different ladders, or equivalently by Lemma~\ref{lem: 321 not sum of ladders}, $\overline{\m}_{x'}$ is not the sum of two ladders: a contradiction. Thus $(i,j-1)\notin Y$.

In the case (B), $b(\Delta_2') = b(\Delta_j^{(\n)})$ and, as $\m_{x'}$ tiles $(\m,\n)$, there is another segment $\Delta_{i'}^{(\m)}$ in $\m$ such that $e(\Delta_{i'}^{(\m)})+1  = b(\Delta_j^{(\n)})$. Then $\Delta_{i'}^{(\m)}\prec_k \Delta_j^{(\n)}$, thus $(i',j)\in X$\footnote{Note that in the illustration above, we have $\Delta_i^{(\m)}\prec_k \Delta_{i'}^{(\m)}$, but it is not necessarily the case, there can be several segments between $\Delta_i^{(\m)}$ and $\Delta_{i'}^{(\m)}$.}. If $\Delta_{j-1}^{(\n)} \in \n$ exists, we have $e(\Delta_{i'}^{(\m)}) > e(\Delta_{i}^{(\m)}) \geq e(\Delta_{j-1}^{(\n)})$, so $(i',j-1)\notin Y$.
In both cases we have a pair $(i,j)\in X$ such that $(i,j-1)\notin Y$ and $\Delta$ tiles the beginning of $\Delta_{i}^{(\m)}$, or a segment before it in $\m+\n$, and tiles the end of $\Delta_j^{(\n)}$.

If $\Delta_{i+1}^{(\m)}$ does not exist, or if $(i+1, j) \notin Y$, then the condition ${\rm NC}_k(\m, \n)$ is satisfied, with data $(i,j,0)$, so ${\rm LC}_k(\m, \n)$ does not hold. 

Suppose $(i+1, j) \in  X$, then we replace $(i,j)$ by $(i+1,j)$. If $\Delta_{j-1}^{(\n)} $ exists, $e(\Delta_{j-1}^{(\n)}) < e(\Delta_{i}^{(\m)})< e(\Delta_{i+1}^{(\m)})$, or $b(\Delta_{j-1}^{(\n)}) \leq b(\Delta_{i}^{(\m)})< b(\Delta_{i+1}^{(\m)})$, so $(i+1, j-1) \notin  Y$. If $\Delta_{i+2}^{(\m)}$ does not exist, or if $(i+2, j) \notin Y$, then the condition ${\rm NC}_k(\m, \n)$ is satisfied, with data $(i+1,j,0)$, so ${\rm LC}_k(\m, \n)$ does not hold. If $(i+2, j) \in X$, replace $(i+1,j)$ by $(i+2,j)$. Recursively, either ${\rm LC}_k(\m, \n)$ does not hold, or there exists $i'\geq i$, such that $(i',j)\in X$, $(i',j-1)\notin Y$, $(i'+1,j) \in Y\setminus X$ and there is segment $\Delta\in \m_{x'}$ which tiles the beginning of $\Delta_{i'}^{(\m)}$, or a segment before it in $\m+\n$, and tiles the end of $\Delta_j^{(\n)}$. Assume the latter. 

\begin{center}
    \begin{tikzpicture}
        \node at (0,0) {\begin{tikzpicture}
\draw[ultra thick] (0,0.5) -- (4,0.5);
\draw[ultra thick] (1.5,1) -- (6,1);
\draw[ultra thick] (1.5,1.5) -- (5,1.5);
\draw[ultra thick, dotted] (1.5,1.5) -- (1.5,1);
\draw[line width=5pt, purple, opacity=0.2] (-0.5,0) -- (3,0);
\draw[line width=5pt, purple, opacity=0.2] (3.2,1) -- (6,1);
\node at (-0.7,0.5) {$\Delta_{i'}^{(\m)}$};
\node at (1,1.5) {$\Delta_{i'+1}^{(\m)}$};
\node at (0,1) {$\Delta_j^{(\n)}$};
\node[purple] at (-1,0) {$\Delta$};
\end{tikzpicture}};
\node at (4.5,0) { or };
 \node at (8,0) {\begin{tikzpicture}
\draw[ultra thick] (0,0.5) -- (4,0.5);
\draw[ultra thick] (1.5,1.5) -- (5,1.5);
\draw[ultra thick] (0.5,1) -- (5,1);
\draw[ultra thick, dotted] (5,1) -- (5,1.5);
\draw[line width=5pt, purple, opacity=0.2] (-0.5,0) -- (3,0);
\draw[line width=5pt, purple, opacity=0.2] (3.2,1.5) -- (5,1.5);
\node at (-0.7,0.5) {$\Delta_{i'}^{(\m)}$};
\node at (0,1) {$\Delta_{i'+1}^{(\m)}$};
\node at (0.8,1.5) {$\Delta_j^{(\n)}$};
\node[purple] at (-1,0) {$\Delta$};
\end{tikzpicture}};
    \end{tikzpicture}
\end{center}

Suppose there is no $\Delta_{j+1}^{(\n)} \in \n$. Since $\m_{x'}$ tiles $(\m,\n)$, there exists $\Delta' \in \m_{x'}$ such that $e(\Delta')=e(\Delta_{i'+1}^{(\m)})$ and $b(\Delta') > b(\Delta)$. By the tiling algorithm, $b(\Delta') \leq b(\Delta_{i'+1}^{(\m)})$. There also exists $\Delta'' \in \m_{x'}$ such that $e(\Delta'')=e(\Delta_{i}^{(\m)})$, and by the tiling algorithm $b(\Delta'') \geq b(\Delta_{i'+1}^{(\m)})$. Hence $b(\Delta) < b(\Delta') \leq b(\Delta_{i'+1}^{(\m)})$ and $e(\Delta) \geq e(\Delta') > e(\Delta_{i}^{(\m)})$, which is a contradiction with the 321-pattern avoidance. So there exists a segment $\Delta_{j+1}^{(\n)} \in \n$. Let $\Delta_1\in \m_{x'}$ be the segment which tiles the end of $\Delta_{j+1}^{(\n)}$. Either $\Delta_1=\Delta$, or $\Delta_1$ is the first segment in $\m_{x'}$ after $\Delta$ with a non-trivial tiling sequence. In any case, there is a segment in $\m_{x'}$ which tiles the beginning of $\Delta_{i'+1}^{(\m)}$, or a segment before it in $\m+\n$, and tiles the end of $\Delta_{j+1}^{(\n)}$.

If $(i'+1,j+1)\in X$, and either $\Delta_{i'+2}^{(\m)}$ does not exist or $(i'+2,j+1)\notin Y$ then the condition ${\rm LC}_k(\m, \n)$ does not hold. As above, if $(i'+2,j+1)\in X$, replace $(i'+1,j+1)$ by $(i'',j+1)$ until $(i'',j+1)\in X$, $(i'',j)\notin Y$, $(i''+1,j+1)\in Y\setminus X$ and there is a segment in $\m_{x'}$ which tiles the beginning of $\Delta_{i''}^{(\m)}$, or a segment before it in $\m+\n$, and tiles the end of $\Delta_{j+1}^{(\n)}$. As above, this implies the existence of $\Delta_{j+2}^{(\n)} \in \n$, and of a segment $\Delta_2\in \m_{x'}$ (possibly equal to $\Delta_1$), which tiles the beginning of $\Delta_{i''+1}^{(\m)}$, or a segment before it in $\m+\n$, and tiles the end of $\Delta_{j+2}^{(\n)}$.

By replacing recursively the segments, either the condition ${\rm LC}_k(\m, \n)$ does not hold or there exists ${\bf i}, {\bf j}$ and $r\geq 0$, such that $({\bf i}, {\bf j}-1) \notin Y$, for all $0\leq l \leq r-1$, $({\bf i} +l, {\bf j}+l) \in X$,  $({\bf i} +l +1 , {\bf j} +l) \in Y \setminus X$, there are segments $\Delta_l\in \m_{x'}$ (not necessarily distinct) which tile the beginning of $\Delta_{{\bf i}+l}^{(\m)}$, or a segment before it in $\m+\n$, and tiles the end of $\Delta_{{\bf j}+l}^{(\n)}$, and a segment $\Delta_r\in \m_{x'}$ which tile the beginning of $\Delta_{{\bf i}+r}^{(\m)}$, or a segment before it in $\m+\n$, and tiles the end of $\Delta_{{\bf j}+r}^{(\n)}$, and $({\bf i} +r  , {\bf j} +r) \notin X$.

Then, necessarily, $e(\Delta_{{\bf i}+r}^{(\m)}) < b(\Delta_{{\bf j}+r}^{(\n)})-1$. Since $\Delta_r$ tiles, there exists $\Delta_{{\bf i'}}^{(\m)}\in \m$ such that $b(\Delta_{{\bf j}+r}^{(\n)})-1 \in \Delta_{{\bf i'}}^{(\m)}$, and $\Delta_r$ tiles part of $\Delta_{{\bf i'}}^{(\m)}$ (similarly to case $(B)$ in the figure above). 
Here we distinguish two cases. If the end of $\Delta_{{\bf i'}}^{(\m)}$ has not been tiled by a previous segment in $\m_{x'}$, then $e(\Delta_{{\bf i'}}^{(\m)}) = b(\Delta_{{\bf j}+r}^{(\n)})-1$, so $(\mathbf{i'},\mathbf{j}+r)\in X$, with $(\mathbf{i'},\mathbf{j}+r-1)\notin Y$, and we can restart the reasoning with the pair $(\mathbf{i'},\mathbf{j}+r)$ (the segment $\Delta_r$ tiles the beginning of a segment before $\Delta_{{\bf i'}}^{(\m)}$ and the end of $\Delta_{{\bf j}+r}^{(\n)}$). 

If the end of $\Delta_{{\bf i'}}^{(\m)}$ has already been tiled by a previous segment in $\m_{x'}$, since we have only been considering segments in $\m_{x'}$ with trivial tiling sequences and the $\Delta_l$'s, necessarily $\Delta_{r-1}$ tiles the end of $\Delta_{{\bf i'}}^{(\m)}$ and ends at $e(\Delta_{{\bf i'}}^{(\m)})$. 

\begin{center} 
\begin{tikzpicture}
\draw[ultra thick] (0,0.5) -- (6,0.5);
\draw[ultra thick] (0,1) -- (3,1);
\draw[ultra thick] (1,1.5) -- (4,1.5);
\draw[ultra thick] (2,2) -- (5,2);
\draw[ultra thick] (4,2.5) -- (9,2.5);
\draw[ultra thick] (6.5,3) -- (10,3);
\draw[ultra thick, dotted] (0,1) -- (0,0.5);
\draw[line width=5pt, purple, opacity=0.2] (-0.5,0) -- (10,0);
\draw[line width=5pt, purple, opacity=0.2] (1,1) -- (3,1);
\draw[line width=5pt, purple, opacity=0.2] (3.5,1.5) -- (4,1.5);
\draw[line width=5pt, purple, opacity=0.2] (4.5,2) -- (5,2);
\draw[line width=5pt, purple, opacity=0.2] (5.5,2.5) -- (6,2.5);
\draw[line width=5pt, purple, opacity=0.2] (6.5,3) -- (10,3);
\draw[line width=5pt, blue, opacity=0.3] (-0.5,-0.5) -- (9,-0.5);
\draw[line width=5pt, blue, opacity=0.3] (2,0.5) -- (6,0.5);
\draw[line width=5pt, blue, opacity=0.3] (6.5,2.5) -- (9,2.5);
\node at (-0.7,0.4) {$\Delta_{{\bf j}+r-1}^{(\n)}$};
\node at (-0.7,1.2) {$\Delta_{{\bf i}+r}^{(\m)}$};
\node at (3.5,2.5) {$\Delta_{{\bf i'}}^{(\m)}$};
\node at (6,3) {$\Delta_{{\bf j}+r}^{(\n)}$};
\node[purple] at (10.5,0) {$\Delta_r$};
\node[blue] at (9.5,-0.5) {$\Delta_{r-1}$};
\end{tikzpicture}   
\end{center}

This implies that $\Delta_{{\bf i'}}^{(\m)}$ is tiled directly after $\Delta_{{\bf j}+r-1}^{(\n)}$ by $\Delta_{r-1}$, and so that $(\mathbf{j} +r, \mathbf{i'})\in X^{(k)}_{\n,\m}$. Moreover, we have $e(\Delta_{{\bf j}+r-1}^{(\n)}) = b(\Delta_{{\bf j}+r}^{(\n)}) -1 $. So if $(\mathbf{j} +r, \mathbf{i'}-1)\in Y^{(k)}_{\n,\m}$, $e(\Delta_{{\bf i'}-1}^{(\m)}) \geq e(\Delta_{{\bf j}+r-1}^{(\n)}) = b(\Delta_{{\bf j}+r}^{(\n)}) -1$, which is a contradiction with the definition of $\mathbf{i'}$. Thus we can restart the reasoning with the pair $(\mathbf{j} +r, \mathbf{i'})\in X^{(k)}_{\n,\m}$, knowing that there is a segment tiling the beginning of $\Delta_{{\bf j}+r}^{(\n)}$, or a segment before it in $\m+\n$, and the end of $\Delta_{{\bf i'}}^{(\m)}$. 

Pursuing recursively this reasoning, we arrive at a pair satisfying either the condition $\mathrm{NC}_k(\m,\n)$, or $\mathrm{NC}_k(\n,\m)$.

\end{proof}

We give the following examples to explain Proposition \ref{prop:ladder times ladder}.

\begin{example}
Let $\m=[-6,-1] + [-2,3] + [-1,4]$, $\n = [-4,1] + [0,2]$ and let $k=15$. We have that $\lambda = (-6,-4,-2,-1,0)$, $\mu=(-1,1,2,3,4)$. Since the stabilizers $S_{\lambda}$ and $S_{\mu}$ are both trivial, we have that for each $x \in S_m$, the only one element $x' \in S_m$ such that $\m_{x'} = \m_x$ is $x' = x$. 

For $x = 31245 = s_1s_2 \in S_5$, $\m_x = [-6, 2] + [-4, -1] + [-2, 1] + [-1, 3] + [0, 4]$. Here $x= 31245$ is the one-line notation of a permutation. The permutation $x$ is $321$-avoiding and $\om_x$ is not $0$. We have
\begin{align*}
{\rm Seq}(\m, \n,  [-6, 2]) = ( ([-6,-1], [-6,-1], \m), ([0,1], [-4,1], \n), ([2,2], [-2,3], \m) ).
\end{align*}
Sine the end points of $[2,2], [-2,3]$ are different, the segment $[-6,2]$ does not tile $\m, \n$. Therefore $\m_x$ does not tile $\m, \n$. Hence $[\oZ(\m_x)]$ does not appear on the right hand side of (\ref{eq:Zm times Zn ladders in Rk}).

For $x = 24153 = s_4s_2s_3s_1 \in S_5$, $\m_x = [-6, 1]+[-4, 3]+[-2, -1]+[-1, 4]+[0, 2]$. The permutation $x$ is $321$-avoiding and $\om_x$ is not $0$. We have
\begin{align*}
& {\rm Seq}(\m, \n,  [-6, 1]) = ( ([-6,-1], [-6,-1], \m), ([0,1], [-4,1], \n) ), \\
& {\rm Seq}(\m_1, \n_1,  [-4,3]) = ( ([-4,-1], [-4,-1], \n), ([0,3], [-2,3], \m) ), \\
& {\rm Seq}(\m_2, \n_2,  [-2, -1]) = ( ([-2,-1], [-2,-1], \m) ), \\
& {\rm Seq}(\m_3, \n_3,  [-1, 4]) = ( ([-1,4], [-1,4], \m) ), \\
& {\rm Seq}(\m_4, \n_4,  [0, 2]) = ( ([0,2], [0,-2], \n).
\end{align*}
Therefore $\m_x$ tiles $\m, \n$. Hence $[\oZ(\m_x)]$ appears on the right hand side of (\ref{eq:Zm times Zn ladders in Rk}).

By checking all $x \in S_5$, we obtain that
\begin{align*}
[\oZ(\m)][\oZ(\n)] & = \oZ( [[-6, -1] + [-4, 1] + [-2, 3] + [-1, 4] + [0, 2]])] \\
& +  \oZ( [[-6, -1] + [-4, 3] + [-2, 1] + [-1, 4] + [0, 2]])] \\ 
& + \oZ( [[-6, 1] + [-4, -1] + [-2, 3] + [-1, 4] + [0, 2]])] \\
& + \oZ( [[-6, 3] + [-4, -1] + [-2, 1] + [-1, 4] + [0, 2]])] \\
& + \oZ( [[-6, 1] + [-4, 3] + [-2, -1] + [-1, 4] + [0, 2]])].
\end{align*}

Now we check the condition ${\rm LC}_k(\m, \n)$ and ${\rm LC}_k(\n, \m)$. We have that 
\begin{align*}
X^{(k)}_{\m,\n}=\{(1, 1), (1, 2)\}, \quad Y^{(k)}_{\m,\n} = \{(1, 1)\}.
\end{align*}
The size of $Y^{(k)}_{\m,\n}$ is less than the size of $X^{(k)}_{\m,\n}$. Therefore there is no injective matching function from $X^{(k)}_{\m,\n}$ to $Y^{(k)}_{\m,\n}$. The condition ${\rm LC}_k(\m, \n)$ is not satisfied.

We have that 
\begin{align*}
X^{(k)}_{\n,\m}=\{(1, 2), (1, 3)\}, \quad Y^{(k)}_{\n,\m} = \{(1, 2), (1,3)\}.
\end{align*}
So the condition ${\rm NC}_k(\n, \m)$ is satisfied, with $(i,j,r)=(1,2,0)$, hence the condition ${\rm LC}_k(\n, \m)$ is not satisfied.

In the case of $k=7$, $[\oZ(\m)][\oZ(\n)] = [\oZ(\m+\n)]$. Indeed, all 
\[X^{(k)}_{\m,\n}=Y^{(k)}_{\m,\n} = X^{(k)}_{\n,\m}=Y^{(k)}_{\n,\m} =  \emptyset,\]
so the conditions ${\rm LC}_k(\m, \n)$ and ${\rm LC}_k(\n, \m)$ are satisfied.
\end{example}

\begin{example}
Let $\m=[-9, -4] + [-8, -2] + [-5, 0]$, $\n = [-10, -3] + [-7, -2] + [-6, -1]$ and let $k \ge 9$.  
Then 
\begin{align*}
[\oZ(\m)][\oZ(\n)] = & \oZ[([-5, 0] + [-6, -1] + [-7, -2] + [-8, -2] + [-9, -4] + [-10, -3])] \\
& +
\oZ[([-5, 0]+ [-6, -1]+ [-7, -2]+ [-8, -3]+ [-9, -4]+ [-10, -2])] \\
& +
\oZ[([-5, -1]+ [-6, 0]+ [-7, -2]+ [-8, -2]+ [-9, -4]+ [-10, -3])] \\
& +
\oZ[([-5, 0]+ [-6, -2]+ [-7, -4]+ [-8, -1]+ [-9, -2]+ [-10, -3])] \\
& + 
\oZ[([-5, -1]+ [-6, 0]+ [-7, -2]+ [-8, -3]+ [-9, -4]+ [-10, -2])] \\
& +
\oZ[([-5, -2]+ [-6, 0]+ [-7, -4]+ [-8, -1]+ [-9, -2]+ [-10, -3])] \\
& +
\oZ[([-5, -1]+ [-6, -2]+ [-7, -4]+ [-8, 0]+ [-9, -2]+ [-10, -3])]. 
\end{align*}
For each $k \ge 9$, $\oZ(\m) \times \oZ(\n)$ is not irreducible. Now we check the conditions ${\rm LC}_k(\m, \n)$ and ${\rm LC}_k(\n, \m)$. We have that 
\begin{align*}
X^{(k)}_{\m,\n}= \{ (1,2), (1,3), (2,3) \}, \ Y^{(k)}_{\m,\n} = \{ (1,2), (1,3), (2,2), (2,3) \}.
\end{align*} 

So $(2,2)\in Y^{(k)}_{\m,\n}\setminus X^{(k)}_{\m,\n}$ and the condition ${\rm NC}_k(\m, \n)$ is satisfied, with $(i,j,r)=(1,2,1)$, hence the condition ${\rm LC}_k(\m, \n)$ is not satisfied. On the other hand:
\begin{align*}
X^{(k)}_{\n,\m}= \{ (1,2), (1,3), (2,3), (3,3)  \}, \ Y^{(k)}_{\n,\m} = \{  (1,2), (1,3), (2,3), (3,3)\}.
\end{align*} 
The condition ${\rm NC}_k(\n, \m)$ is satisfied, with $(i,j,r)=(1,2,0)$, hence the condition ${\rm LC}_k(\n, \m)$ is not satisfied.


\end{example}

\section{Application to Grassmannian cluster algebras} \label{sec:application to Grassmannian cluster algebras}

In this section, we apply our criterion (Theorem \ref{thm:condition of irreducibility of tensor product}) of simplicity of tensor product of two snake modules to obtain compatibility criterion of two cluster variables in a Grassmannian cluster algebra which correspond to snake modules. This result generalizes Leclerc and Zelevinsky's result about compatibility of two Pl\"{u}cker coordinates \cite{LZ}. 

\subsection{Grassmannian cluster algebras and semistandard Young tableaux}  \label{subsec:Grassmannian cluster algebras and semistandard Young tableaux}
Let us recall some facts about semi-standard Young tableaux and their relation to the Grassmannian cluster algebra. 


A {\sl semistandard Young tableau} is a Young tableau with weakly increasing rows and strictly increasing columns. For $k,n \in \ZZ_{\ge 1}$, we denote by ${\rm SSYT}(k, [n])$ the set of rectangular semistandard Young tableaux with $k$ rows and with entries in $[n]$ (with arbitrarly many columns). The empty tableau is denoted by $\mathds{1}$. 

For $S,T \in {\rm SSYT}(k, [n])$, let $S \cup T$ be the row-increasing tableau whose $i$th row is the union of the $i$th rows of $S$ and $T$ (as multisets), \cite{CDFL}.  

We call $S$ a factor of $T$, and write $S \subset T$, if the $i$th row of $S$ is contained in that of $T$ (as multisets), for $i \in [k]$. In this case, we define $\frac{T}{S}=S^{-1}T=TS^{-1}$ to be the row-increasing tableau whose $i$th row is obtained by removing that of of $S$ from that of $T$ (as multisets), for $i \in [k]$. 

A tableau $T \in {\rm SSYT}(k, [n])$ is \emph{trivial} if each entry of $T$ is one less than the entry below it. 

For any $T \in {\rm SSYT}(k, [n])$, we  denote by $T_{\text{red}} \subset T$ the semistandard tableau obtained by removing a maximal trivial factor from $T$. For trivial $T$, one has $T_{\text{red}} = \mathds{1}$. 

Let ``$\sim$'' be the equivalence relation on $S, T \in {\rm SSYT}(k, [n])$ defined by: $S \sim T$ if and only if $S_{\text{red}} = T_{\text{red}}$. We denote by ${\rm SSYT}(k, [n],\sim)$ the set of $\sim$-equivalence classes.

Denote by $\Gr(k,n)$ the Grassmannian of $k$-planes in $\mathbb{C}^n$ and $\CC[\Gr(k,n)]$ its homogeneous coordinate ring. Define $\CC[\Gr(k,n,\sim)]$ to be the quotient of $\CC[\Gr(k,n)]$ by the ideal generated by $P_{i_1, \ldots, i_n}-1$ with $\{i_1, \ldots, i_n\}$ being a consecutive interval, and $P_{i_1, \ldots, i_n} \in \CC[{\rm Gr}(k, n)]$ being the Pl\"ucker coordinate.  

In \cite{CDFL}, it is shown that the elements in the dual canonical basis of $\CC[\Gr(k,n,\sim)]$ are in bijection with semistandard Young tableaux in ${\rm SSYT}(k, [n],\sim)$. The elements in the dual canonical basis of $\CC[\Gr(k,n,\sim)]$ are in bijection with simple modules in a certain category of finite dimensional $U_q(\widehat{\mathfrak{sl}_k})$-modules \cite{HL10}.

A one-column tableau is called a \emph{fundamental tableau} if its entries are a $k$-set of the form $[i,i+k] \setminus \{r\}$ for $r \in \{i+1, \ldots, i+k-1\}$. A tableau $T$ is said to have \emph{small gaps} if each of its columns is a fundamental tableau. Then any tableau in $\SSYT(k,[n])$ is $\sim$-equivalent to a unique small gap semistandard tableau.


We now recall the explicit formula of $\ch(T)$ \cite{CDFL}. For that we need to first define $w_T \in S_m$, $P_{u; T'}$, $u \in S_m$, for every $T \in \SSYT(k, [n])$, where $T'$ is the unique small gap tableau which is $\sim$-equivalent to $T$, and $m$ is the number of columns of $T'$. 

Let ${\bf i} = i_1 \leq i_2 \dots \leq i_m$ be the entries in the first row of $T'$, and let $r_1,\dots,r_m$ be the elements such that the $a$th column of $T'$ has content $[i_a,i_a+n] \setminus \{r_a\}$. Let ${\bf j} = j_1 \leq j_2 \leq \dots \leq j_m$ be the elements $r_1,\dots,r_m$ written in weakly increasing order. 

For $u \in S_m$, define $P_{u;T'} \in \CC[\Gr(k,n)]$ as follows. Provided $j_a \in [i_{u(a)}, i_{u(a)}+k]$ for all $a \in [m]$, define the tableau $\alpha(u;T')$ to be the semistandard tableau whose columns have content $[i_{u(a)}, i_{u(a)}+k] \setminus \{j_a\}$ for $a \in [m]$, and define $P_{u; T'} = P_{\alpha(u;T')} \in \CC[\Gr(k,n)]$ to be the corresponding standard monomial. On the other hand, if $j_a \notin [i_{u_a}, i_{u(a)}+k]$ for some $a$, then the tableau $\alpha(u;T')$ is {\sl undefined} and $P_{u ;T'} = 0$. 

There is a unique $u \in S_m$ is of maximal length with the property that the sets $\{[i_{u(a)},i_{u(a)}+k] \setminus \{j_a \} \}_{a \in [m]}$ describe the columns of $T'$. This $u$ is denoted by $u = w_{T}$.

By \cite[Theorem 5.8]{CDFL}, the element $\ch(T)$ in the dual canonical basis of $\CC[\Gr(k,n,\sim)]$ is given by 
\begin{align}\label{eq:formula of ch(T)}
\ch(T) = \sum_{u \in S_m} (-1)^{\ell(uw_T)} p_{uw_0, w_Tw_0}(1) P_{u; T'} \in \CC[\Gr(k,n,\sim)]
\end{align}
where $p_{u,v}(q)$ is a Kazhdan-Lusztig polynomial.

\subsection{Tableaux and multisegments} \label{subsec:tableaux and multisegments}

For $a,b\in\ZZ$ with $a\leq b < a +k-1$, denote by $T_{a,b}$ the fundamental tableau with entries $\{1-a,2-a, \ldots, k-a+1\} \backslash \{k-b\}$. Denote by ${\rm Fund}(k, [n])$ the set of fundamental tableaux with $k$ rows and with entries in $[n]$. From \cite{CDFL}, the tableau $T_{a,b}$ corresponds to the segment $[a,b]$. We have:
\begin{equation}\label{eq:multisegtoTableaux}
\begin{array}{ccc}
\mathrm{Seg}_k :=\{ [a,b]\in\mathrm{Seg}_k \mid k+1-n\leq a\leq 0\}  & \simeq & {\rm Fund}(k, [n]),\\
 {[}a, b] & \mapsto & T_{a,b},\\
 {[}1-i, k-r] & \mapsfrom & T_{ \{i, i+1 ,\ldots,i+k-1\} \setminus \{r\}}. 
\end{array}
\end{equation}
Any tableau $T$ can be decomposed uniquely into the union of fundamental tableaux. We define the multisegment ${\bf m}_T$ corresponding to $T$ as the sum of the segments corresponding to the fundamental tableaux in the decomposition. For a multisegment ${\bf m}$, we define $T_{\bf m}$ to be the union of the fundamental tableaux corresponding to the segments in ${\bf m}$.


By the correspondence between dominant monomials and multisegments in \eqref{eq:multisegtomonom}, we have a correspondence between dominant monomials and tableaux induced by the following correspondence:
\begin{equation}\label{eq:monomToTableaux}
\begin{array}{ccc}
\{ Y_{i,p} \mid (i,p)\in \hat{I}, 1\leq i+p\leq 2n-2k \} & \simeq & {\rm Fund}(k, [n]),\\
 Y_{i,p} & \mapsto & T_{\frac{1-i-p}{2}, \frac{i-p-1}{2}}, \\
 Y_{k-r+i, r-k+i-1} & \mapsfrom & T_{ \{i, i+1 ,\ldots,i+k-1\} \setminus \{r\}}.
\end{array}
\end{equation}

We denote the tableau corresponding to a dominant monomial $M$ by $T_M$ and denote by $M_T$ the dominant monomial corresponding to a tableau $T$.  

For a one-column tableau $T$ with entries $i_1, \ldots, i_k$, we call the numbers in $[i_1, i_k] \setminus \{i_1, \ldots, i_k\}$ the \emph{missing numbers}. 

For a one-column tableau $T$, we define two sequences ${\bf i}_T$, ${\bf j}_T$ as follows. If the set of missing numbers of $T$ is empty, then ${\bf i}_T$, ${\bf j}_T$ are empty. If the set of missing numbers of $T$ is not empty, then ${\bf j}_T$ is the increasing sequence with entries being the missing numbers of $T$ and ${\bf i}_T$ is $(i_1, \ldots, i_1+r-1)$, where $i_1$ is the first entry of $T$ and $r$ is the size of the set of missing numbers of $T$. 

For a tableau $T$ with columns $T_1, \ldots, T_m$, we define ${\bf i}_T = {\bf i}_{T_1} \cdots {\bf i}_{T_m}$, ${\bf j}_T = {\bf j}_{T_1} \cdots {\bf j}_{T_m}$ to be the concatenation of sequences. 

A multisegment $\m=\sum_{i=1}^m [a_i, b_i]$ is called \emph{regular} if $a_i\neq a_j$ and $b_i\neq b_j$, for all $i\neq j$. A multisegment $\m=\sum_{i=1}^m [a_i, b_i]$ is said to be a \emph{ladder} if $a_1< \ldots < a_m$ and $b_1 < \ldots < b_m$. 

The following is clear.

\begin{lemma}
For a tableau $T$, ${\bf m}_T$ is regular if and only if the numbers in ${\bf i}_T$ are all different and the numbers in ${\bf j}_T$ are all different. 

The multisegment ${\bf m}_T$ is a ladder if and only if the numbers in ${\bf j}_{T_r}$ are all less than the numbers in ${\bf j}_{T_{r+1}}$ for every $r=1,\ldots,m-1$, where $T_1, \ldots, T_m$ are columns of $T$.
\end{lemma}

When ${\bf m}_T$ is a ladder, we call $T$ and $\ch(T)$ a ladder or a snake. The $U_q(\widehat{\mathfrak{sl}_k})$-module $L(M_{\bf m})$ corresponding to ${\bf m}_T$ is a \emph{snake module}, as introduced in \cite{MY12}.

\subsection{Weakly separated and unlinked property}
\begin{definition}[{\cite{LZ}}] \label{def:weakly separted k-subsets}
Given two $k$-subsets $I$ and $J$ of $\{1, \dots, n\}$, denote by $\min(J)$ the minimal element in $J$ and by $\max(I)$ the maximal element in $I$, we write $I < J$ if $\max(I)<\min(J)$. The sets $I$ and $J$ are called \emph{weakly separated} if at least one of the following two conditions holds:
\begin{enumerate}
\item $J\setminus I$ can be partitioned into a disjoint union $J\setminus I = J' \sqcup J''$ so that $J' < I\setminus J < J''$;

\item $I\setminus J$ can be partitioned into a disjoint union $I\setminus J = I' \sqcup I''$ so that $I' < J\setminus I < I''$.
\end{enumerate}
\end{definition}

Recall that in Section \ref{subsec:tableaux and multisegments}, we denote $T_{a,b}$ to be the one-column tableau with entries $1-a, 2-a, \ldots, k-b-1, k-b+1, \ldots, k-a+1$. Two fundamental tableaux are called \emph{weakly separated} (w.s.) if the corresponding $k$-subsets are weakly separated.

We have the following easy criterion to determine if two one column tableaux with small gaps are weakly separated. 
\begin{lemma}\label{lem:criterion_n_w_s}
Let $T,T'$ be one column tableaux with small gaps. Write $T=T_{a,b}$ and $T'=T_{c,d}$, and let $I_T,I_{T'}$ be the corresponding $k$-subsets of $[n]$. Then $T$ and $T'$ are \emph{not weakly separated} if and only if
\[\left\lbrace \begin{array}{l}
k-b\in I_{T'},\\
k-d \in I_T,
\end{array}\right. \quad  \text{and} \quad (a-c)(b-d)>0.\]
\end{lemma}

\begin{proof}
For the direct implication, suppose without loss of generality that $a\leq c$. Then let us suppose first that $b=d$. If $a=c$ then $T=T'$ and $T,T'$ are weakly separated. Thus $a<c$, then
\begin{align*}
& I_T\setminus I_{T'} = \{ k-c+2 , k-c+3,\ldots,k-a+1\},\\
\text{and } & I_{T'}\setminus I_T = \{1-c,2-c,\ldots, -a\}.
\end{align*}
We have $I_{T'}\setminus I_T < I_T\setminus I_{T'}$. Thus $T,T'$ are w.s.

Now suppose that $k-b \notin I_{T'} = \{ 1-c, 2-c,\ldots , k-c+1\}$. As $k-b > 1-a \geq 1-c$, we have $k-b > k-c+1$. We consider two subcases: $k-d< 1-a$ and $k-d\geq 1-a$, as illustrated below.

\begin{center}
\begin{tikzpicture}
\draw[|-|] (0,0) to (1.9,0);
\draw[|-|] (2.1,0) to (3,0);
\draw[|-|] (-2,-0.3) to (-1.1,-0.3); 
\draw[|-|] (-0.9,-0.3) to (1,-0.3); 
\draw[|-|] (1.1,-0.8) to (1.9,-0.8); 
\draw[|-|] (2.1,-0.8) to (3,-0.8); 
\draw[|-|] (-2,-1.1) to (-1.1,-1.1); 
\draw[|-|] (-0.9,-1.1) to (-0.1,-1.1);
\node at (-2.7,0.1) {\small $I_T$};
\node at (-2.7,-0.3) {\small $I_{T'}$}; 
\node at (0,0.2) {\tiny $1-a$};
\node at (2,0.2) {\tiny $k-b$};
\node at (3.2,0.2) {\tiny $k-a+1$};
\node at (-2,-0.1) {\tiny $1-c$};
\node at (-1,-0.5) {\tiny $k-d$};
\node at (1,-0.5) {\tiny $k-c+1$};
\node at (-2.7,-0.7) {\small $I_T\setminus I_{T'}$};
\node at (-2.7,-1.1) {\small $I_{T'}\setminus I_T$};
\end{tikzpicture}
\begin{tikzpicture}
\draw[|-|] (-1,0) to (1.9,0);
\draw[|-|] (2.1,0) to (3,0);
\draw[|-|] (-2,-0.3) to (-0.1,-0.3); 
\draw[|-|] (0.1,-0.3) to (1,-0.3); 
\draw[|-|] (1.1,-0.8) to (1.9,-0.8); 
\draw[|-|] (2.1,-0.8) to (3,-0.8); 
\draw[|-|] (-2,-1.1) to (-1.1,-1.1); 
\node at (0,-0.8) {\tiny $\bullet$};
\node at (-2.7,0.1) {\small $I_T$};
\node at (-2.7,-0.3) {\small $I_{T'}$}; 
\node at (-1,0.2) {\tiny $1-a$};
\node at (2,0.2) {\tiny $k-b$};
\node at (3.2,0.2) {\tiny $k-a+1$};
\node at (-2,-0.1) {\tiny $1-c$};
\node at (0,-0.5) {\tiny $k-d$};
\node at (1.2,-0.5) {\tiny $k-c+1$};
\node at (-2.7,-0.7) {\small $I_T\setminus I_{T'}$};
\node at (-2.7,-1.1) {\small $I_{T'}\setminus I_T$};
\end{tikzpicture}
\end{center}
If $k-d<1-a$, then 
\begin{align*}
& I_{T} \setminus I_{T'} = [k-c+2, k-b-1] \cup [k-b+1, k-a+1], \\
& I_{T'} \setminus I_T = [1-c, k-d-1] \cup [k-d+1, -a].
\end{align*}
If $k-d>1-a$, then 
\begin{align*}
& I_{T} \setminus I_{T'} = \{k-d\} \cup [k-c+2, k-b-1] \cup [k-b+1, k-a+1], \\
& I_{T'} \setminus I_T = [1-c, -a].
\end{align*}
In both cases, $I_{T'}\setminus I_T < I_T\setminus I_{T'}$, therefore $T$ and $T'$ are w.s.

By a similar reasoning, we can prove that if $k-d\notin \{1-a,2-a,\ldots, k-a+1\}$, then $T$ and $T'$ are w.s.

For the reverse implication, suppose $k-b\in I_{T'}$, $k-d\in I_T$ and $(a-c)(b-d)>0$. We can assume that $a< c$. We have $k-b\leq k-c+1$ and then the configuration is as following. 
\begin{center}
\begin{tikzpicture}
\draw[|-|] (-1,0) to (0.9,0);
\draw[|-|] (1.1,0) to (3,0);
\draw[|-|] (-2,-0.3) to (-0.1,-0.3); 
\draw[|-|] (0.1,-0.3) to (2,-0.3); 
\draw[|-|] (2.1,-0.9) to (3,-0.9); 
\draw[|-|] (-2,-1.2) to (-1.1,-1.2); 
\node at (0,-0.9) {\tiny $\bullet$};
\node at (1,-1.2) {\tiny $\bullet$};
\node at (-2.7,0.1) {\small $I_T$};
\node at (-2.7,-0.3) {\small $I_{T'}$}; 
\node at (-1,0.2) {\tiny $1-a$};
\node at (1,0.2) {\tiny $k-b$};
\node at (3,0.2) {\tiny $k-a+1$};
\node at (-2,-0.1) {\tiny $1-c$};
\node at (0,-0.5) {\tiny $k-d$};
\node at (2,-0.5) {\tiny $k-c+1$};
\node at (-2.7,-0.8) {\small $I_T\setminus I_{T'}$};
\node at (-2.7,-1.2) {\small $I_{T'}\setminus I_T$};
\end{tikzpicture}
\end{center}
We have that  
\begin{align*}
& I_{T} \setminus I_{T'} = \{k-d\} \cup [k-c+2, k-a+1], \\
& I_{T'} \setminus I_T = [1-c, -a] \cup \{k-b\}.
\end{align*}
Then $k-d \in I_{T} \setminus I_{T'}$ and there are $x,y \in I_{T'} \setminus I_T$ such that $x< k-d < y$. Similarly, $k-b \in I_{T'} \setminus I_{T}$ is in between elements of $I_{T} \setminus I_{T'}$.
Therefore $T$ and $T'$ are not weakly separated.
\end{proof}

The following result gives an explicit correspondence between two one column tableaux with small gaps not being weakly separated, and their corresponding segment being linked.
\begin{proposition} \label{prop:weakly separated and non linked property}
Let $T,T'$ be one column tableaux with small gaps. Write $T=T_{a,b}$ and $T=T'_{c,d}$, and let $I_T,I_{T'}$ be the corresponding $k$-subsets of $[n]$.  Then the $k$-subsets $I_T$ and $I_{T'}$ are not weakly separated if and only if the corresponding segments $[a,b]$ and $[c,d]$ are linked and 
\[k>\max\left(d-a,b-c\right).\]
\end{proposition}

\begin{proof}
Suppose without loss of generality that $a \le c$. 

Let us suppose that $[a,b]\prec [c,d]$ and $k>d-a$. As $a+1\le c \le b+1 \le d$, then $1-c, k-b\in J\setminus I$, $k-d,k-a+1\in I\setminus J$, and
\[1-c < k-d < k-b < k-a+1.\]
Therefore $I$ and $J$ are not weakly separated. 

Suppose now that $I_T$ and $I_{T'}$ are not weakly separated. From Lemma~\ref{lem:criterion_n_w_s}, $k-d\in [1-a,k-a+1]\setminus\{k-b\}$, thus $k>d-a$. 
Similarly, $k-b\in [1-c,k-c+1]\setminus\{k-d\}$ implies $c\leq b+1$. Additionally, as $(a-c)(b-d)>0$, we deduce that in fact $a<c$ and $d<d$. Therefore $a+1 \le c \le b+1 \le d$, and $[a,b]$ and $[c,d]$ are linked.
\end{proof}

\begin{remark}
We need the condition that $k$ is sufficient large in Proposition \ref{prop:weakly separated and non linked property}. For example, in the case of $k=6$, the multisegments $[-1,2]+[-4,-1]$ correspond to the $k$-subsets $\{2,3,5,6,7,8\}$, $\{5,6,8,9,10,11\}$. The multisegments are linked but the $k$-subsets are weakly separated. 

In the case of $k=7$, we have that the multisegments $[-1,2]+[-4,-1]$ correspond to the $k$-subsets $\{2,3,4,6,7,8,9\}$, $\{5,6,7,9,10,11,12\}$. The multisegments are linked and the $k$-subsets are not weakly separated.
\end{remark}

\subsection{The combinatorial criterion \texorpdfstring{${\rm LC}(T, T')$}{LC(T,T')}} \label{subsec:condition LCTTprime}

Denote by $<$ the lexicographical order on the set of one-column tableaux with gap weight $\le 1$. For example, in ${\rm SSYT}(3,[5])$,
\begin{align*}
\scalemath{0.6}{
\begin{ytableau}
1 \\
2 \\
3 
\end{ytableau} < \begin{ytableau}
1 \\
2 \\
4 
\end{ytableau}
< 
\begin{ytableau}
1 \\
3 \\
4 
\end{ytableau}
< 
\begin{ytableau}
2 \\
3 \\
4 
\end{ytableau} <
\begin{ytableau}
2 \\
3 \\
5 
\end{ytableau}<
\begin{ytableau}
2 \\
4 \\
5 
\end{ytableau}
< 
\begin{ytableau}
3 \\
4 \\
5 
\end{ytableau}. }
\end{align*}

\begin{remark}
The lexicographical order on one-column tableaux with small gaps does not correspond to either the left or right aligned order introduced in Section~\ref{sect_p_adic}. However, if the segments in a multisegments $\m$ are ordered such that the corresponding one column tableaux are in increasing lexicographical order, then the multisegment $\m$ is ordered.

In particular, if two segments $\Delta,\Delta'$ are such that $T_\Delta \leq T_{\Delta'}$ and $\Delta \le_b\Delta'$ (or $\Delta \le_e\Delta'$), then $\Delta$ and $\Delta'$ are unlinked.
\end{remark}

Let $T, T' \in {\rm SSYT}(k, [n])$ with decomposition into unions of one-column small gap tableaux $T=T_1 \cup \cdots \cup T_r$, $T'=T'_1 \cup \cdots \cup T'_{r'}$, where the factors are written in the lexicographical order $T_i < T_{i+1}$, $T'_{i}<T'_{i+1}$.

Denote by ${\rm pr}(T)$ the promotion of $T$ (cf. \cite{Sch72}) and
\begin{align*}
& X_{T,T'} = \{(i,j): T_i, T'_j \text{ are not weakly separated}, \min T'_j < \min T_i\}, \\
& Y_{T,T'} = \{(i,j): \pr(T_i), T'_j \text{ are not weakly separated}, \min T'_j \le \min T_i \}.
\end{align*}

We denote by ${\rm LC}(T, T')$ the following condition: there is an injective map $f: X_{T,T'} \to Y_{T,T'}$, $f(i,j)=(i',j')$, such that for any $(i,j) \in X_{T,T'}$, either $i=i'$, $T'_{j'}$ and $T'_{j}$ are not weakly separated, and $\min T'_j < \min T'_{j'}$ or $j=j'$, $T_i$ and $T_{i'}$ are not weakly separated, and $\min T_{i'}< \min T_i$.


\begin{lemma}\label{lemma: XXYY}
For $T, T' \in {\rm SSYT}(k, [n])$, and ${\bf m}={\bf m}_T$, ${\bf m}'={\bf m}_{T'}$ the corresponding multisegments. Then
\begin{align}
X_{T,T'} &= X^{(k)}_{{\bf m}, {\bf m}'},\\
Y_{T,T'} &= Y^{(k)}_{{\bf m}, {\bf m}'}.
\end{align}
\end{lemma}

\begin{proof}
Let us write the decompositions into one-column small gap tableaux $T=T_1 \cup \cdots \cup T_r$, $T'=T'_1 \cup \cdots \cup T'_{r'}$ and 
let $\Delta_1,\Delta_2,\ldots, \Delta_r$, $\Delta'_1,\Delta'_2,\ldots, \Delta'_{r'}$ be the corresponding segments.

Then the equality $X_{T,T'} = X^{(k)}_{{\bf m}, {\bf m}'}$ is a direct consequence of Proposition~\ref{prop:weakly separated and non linked property}.

Let $(i,j)\in Y_{T,T'}$. If $\pr(T_i)=T_i+1$, where all the values are increased by one, then as before using Proposition~\ref{prop:weakly separated and non linked property} we deduce that $(i,j)\in  Y^{(k)}_{{\bf m}, {\bf m}'}$. 

If the maximal value of $T_i$ is $n$, then $\Delta_i=[a,b]$ with $a=k-n+1$, and the entries of $\pr(T_i)$ are $I=\{1\}\cup[2-a,n]\setminus\{k-b+1\} $. Write $\Delta'_j=[c,d]$. Then the entries of $T_j$ are $J = [1-c, k-c+1]\setminus \{k-d\}$. By assumption $\min(T_j')\leq \min(T_i)$ so $a\leq c$. If $b>d$, then the configuration of segments is the following.

\begin{center}
\begin{tikzpicture}
\node at (-2,0) {\tiny $\bullet$};
\draw[|-|] (0,0) to (0.9,0);
\draw[|-|] (1.1,0) to (4,0);
\draw[|-|] (-1,-0.3) to (1.9,-0.3); 
\draw[|-|] (2.1,-0.3) to (3,-0.3); 
\draw[|-|] (3.1,-0.9) to (4,-0.9); 
\draw[|-|] (-1,-1.2) to (-0.1,-1.2); 
\node at (-2,-0.9) {\tiny $\bullet$};
\node at (2,-0.9) {\tiny $\bullet$};
\node at (1,-1.2) {\tiny $\bullet$};
\node at (-3,0.1) {\small $I$};
\node at (-3,-0.3) {\small $J$};
\node at (-2,0.2) {\tiny $1$}; 
\node at (0,0.2) {\tiny $2-a$};
\node at (1,0.3) {\tiny $k-b+1$};
\node at (4,0.2) {\tiny $n$};
\node at (-1,-0.1) {\tiny $1-c$};
\node at (2,-0.5) {\tiny $k-d$};
\node at (3.2,-0.5) {\tiny $k-c+1$};
\node at (-3,-0.8) {\small $I\setminus J$};
\node at (-3,-1.2) {\small $J\setminus I$};
\end{tikzpicture}
\end{center}
Let $I'=\{1\}$ if $c<0$ and $I'=\emptyset$ otherwise, and $I''= \{k-d\}\cup [k-c+2,n]$ if $c>k-n+1$ and $I''=\{k-d\}$ otherwise. Then 
\begin{align*}
& I \setminus J = I' \sqcup I'', \\
& I' < J \setminus I < I''.
\end{align*}
This contradicts the fact that $\pr(T_i)$ and $T_j'$ are not weakly separated, thus $b\leq d$. Now suppose, $b<c$. Then, as $k-c+1\leq k-b$, we have that $J \setminus I = [1-c, 1-a]$. Therefore all the values in $J\setminus I$ are between 1 and $2-a$ and $\pr(T_i)$ and $T_j'$ are weakly separated: a contradiction. Thus we have $a\leq c\leq b\leq d$, $[a-1,b-1]\prec[c,d]$, and $(i,j)\in Y^{(k)}_{{\bf m}, {\bf m}'}$.

Conversely, let $(i,j)\in Y^{(k)}_{{\bf m}, {\bf m}'}$. If $b(\overleftarrow{\Delta_i}) > k-n+1$, then $\overleftarrow{\Delta_i}$ is the segment corresponding to the tableau $\pr(T_i)=T_i+1$, and we deduce that $(i,j)\in Y_{T,T'}$. Otherwise, the segment $\overleftarrow{\Delta_i}$ has no corresponding tableau. However, if we write $\Delta_i=[a,b]$ and $\Delta'_j=[c,d]$, then the condition $\overleftarrow{\Delta_i}\prec \Delta'_j$ is equivalent to $a\leq c\leq b\leq d$. Thus, the values $k-d$ and $n$ appear in the tableau $\pr(T_i)$ but not in $T'_j$ and the values $1-c$ and $k-b+1$ appear in the tableau $T'_j$ but not in $\pr(T_i)$. As $ 1-c < k-d < k-b+1 < n$, we deduce that $\pr(T_i)$ and $T'_j$ are weakly separated. Thus $(i,j)\in Y_{T,T'}$.
\end{proof}

Using Proposition~\ref{prop:weakly separated and non linked property} and Lemma~\ref{lemma: XXYY}, it is clear that for $T,T'$ tableaux in $\SSYT(k,[n])$ and ${\bf m}_T, {\bf m}_{T'}$ the corresponding multisegments, then
\begin{equation}\label{eq: LC_TT LCmm}
{\rm LC}(T,T') \quad \Leftrightarrow \quad {\rm LC}_k({\bf m}_T, {\bf m}_{T'}).
\end{equation}

\begin{corollary}\label{cor:main result tableaux}
Let $T,T'\in{\rm SSYT}(k,[n])$ be tableaux corresponding to ladders. Then the following conditions are equivalent:
\begin{enumerate}
\item ${\rm LC}(T,T')$ and ${\rm LC}(T', T)$,
\item $\ch(T)\ch(T')=\ch(T\cup T') \in {\rm SSYT}(k,[n],\sim)$.
\end{enumerate} 
\end{corollary}

\begin{proof}
Thanks to \eqref{eq: LC_TT LCmm},$(1)$ is equivalent to both ${\rm LC}_k({\bf m}_T, {\bf m}_{T'})$ and ${\rm LC}_k({\bf m}_{T'}, {\bf m}_T)$. These are equivalent, by Theorem~\ref{thm:condition of irreducibility of tensor product}, to the tensor product $L(M_T)\otimes L(M_{T'})$ being irreducible, or also $\chi_q(L(M_T))\chi_q(L(M_{T'}))=\chi_q(L(M_T M_{T'}))$. The latter is equivalent to $(2)$ by Theorem 3.17 and Proposition 3.26 in \cite{CDFL}. 
\end{proof}

\begin{remark}
We expect Corollary~\ref{cor:main result tableaux} to be true not only in $\CC[\Gr(k,n,\sim)]$ but also in $\CC[\Gr(k,n)]$. In Section 5.2 of \cite{CDFL}, it is conjectured that $\widetilde{\ch}(T) = \frac{1}{P_{T''}} \ch(T)$ is an element in $\CC[\Gr(k,n)]$, where $T'' = T' T^{-1}$, $T'$ is the unique small gap tableau such that $T' \sim T$. 
Thus we expect that condition $(1)$ in Corollary~\ref{cor:main result tableaux} should also be equivalent to 
\[\widetilde{\ch}(T)\widetilde{\ch}(T') = \widetilde{\ch}(T\cup T') \quad \in \CC[\Gr(k,n)].\] 
\end{remark}

\begin{example} Let us look at Example \ref{example:fundamental times fundamental} in the context of tableaux.

    Denote by $T^{(k)}_{{\bf m}}$ the tableau in $\SSYT(k, [n])$ ($n$ is a large enough integer). Since the maximum length of segments in ${\bf m}$, ${\bf m}'$ is $4$, we take $k \ge 5$. 

For $k \ge 5$, ${\rm LC}( T^{(k)}_{{\bf m}'}, T^{(k)}_{{\bf m}})$ is satisfied.
 
For $k=5, 6$, ${\rm LC}( T^{(k)}_{{\bf m}}, T^{(k)}_{{\bf m}'} )$ is satisfied. For $k \ge 7$, ${\rm LC}( T^{(k)}_{{\bf m}}, T^{(k)}_{{\bf m}'} )$ is not satisfied. 

For example, in the case of $k=5$, we have 
\begin{align*}
\scalemath{0.6}{
T^{(k)}_{\bf m} = \begin{ytableau}5\\ 7\\ 8\\ 9\\ 10 \end{ytableau}, \  T^{(k)}_{{\bf m}'} =\begin{ytableau}2\\ 4\\ 5\\ 6\\ 7 \end{ytableau}, }
\end{align*}
and $ X_{ T^{(k)}_{{\bf m}}, T^{(k)}_{{\bf m}'} } = Y_{ T^{(k)}_{{\bf m}}, T^{(k)}_{{\bf m}'}}  = X_{ T^{(k)}_{{\bf m}'}, T^{(k)}_{\bf m} } = Y_{ T^{(k)}_{{\bf m}'}, T^{(k)}_{\bf m} } = \emptyset$. Therefore ${\rm LC}( T^{(k)}_{{\bf m}}, T^{(k)}_{{\bf m}'} )$ and ${\rm LC}( T^{(k)}_{{\bf m}'}, T^{(k)}_{\bf m} )$ are satisfied. 

In the case of $k=7$, we have 
\begin{align*}
\scalemath{0.6}{
T^{(k)}_{\bf m} =\begin{ytableau}5\\ 6\\ 7\\ 9\\ 10\\ 11\\ 12 \end{ytableau}, \  T^{(k)}_{{\bf m}'} = \begin{ytableau}2\\ 3\\ 4\\ 6\\ 7\\ 8\\ 9 \end{ytableau}, }
\end{align*}
and $X_{ T^{(k)}_{{\bf m}}, T^{(k)}_{{\bf m}'} } = \{(1,1)\}$, $Y_{ T^{(k)}_{{\bf m}}, T^{(k)}_{{\bf m}'}} = \emptyset$, $X_{ T^{(k)}_{{\bf m}'}, T^{(k)}_{\bf m} } = Y_{ T^{(k)}_{{\bf m}'}, T^{(k)}_{\bf m} } = \emptyset$. Therefore ${\rm LC}( T^{(k)}_{{\bf m}}, T^{(k)}_{{\bf m}'} )$ is not satisfied and ${\rm LC}( T^{(k)}_{{\bf m}'}, T^{(k)}_{\bf m} )$ is satisfied. 
\end{example}

\begin{example}
    As in Example \ref{example:2nd in main}:
\begin{align*}
& {\bf m} = [-4, -3] + [-5, -4], \\
& {\bf n} = [0, 1] + [-1, 0]+ [-2, -2]  + [-2, -1] + [-3, -3] + [-3, -3] + [-5, -4]. 
\end{align*}

For $k=3$, we have  
\begin{align*}
\ch(T_{\m}^{(3)}) \ch(T_{\n}^{(3)}) = \ch( \scalemath{0.66}{ \begin{ytableau} 5 \\ 8 \\ 9 \end{ytableau} } ) \ch( \scalemath{0.66}{ \begin{ytableau} 1 & 3 & 6\\ 4 & 5 & 8\\ 7 & 7 & 9 \end{ytableau} } ) = \ch( \scalemath{0.66}{ \begin{ytableau} 1 & 3 & 5 & 6\\ 4 & 5 & 8 & 8\\ 7 & 7 & 9 & 9 \end{ytableau} } ) = \ch( T_{\m}^{(3)} \cup T_{\n}^{(3)} ).
\end{align*}  
Now we check the conditions ${\rm LC}(T_{\m}^{(3)}, T_{\n}^{(3)})$ and ${\rm LC}(T_{\n}^{(3)}, T_{\m}^{(3)})$. We have 
\begin{align*}
& \scalemath{0.66}{ T_{\m}^{(3)} \sim  \begin{ytableau} 5\\ 7\\ 8 \end{ytableau}\cup \begin{ytableau} 6\\ 8\\ 9 \end{ytableau}, \quad T_{\n}^{(3)} \sim  \begin{ytableau} 1\\ 3\\ 4 \end{ytableau}\cup \begin{ytableau} 2\\ 4\\ 5 \end{ytableau}\cup \begin{ytableau} 3\\ 4\\ 6 \end{ytableau}\cup \begin{ytableau} 3\\ 5\\ 6 \end{ytableau}\cup \begin{ytableau} 4\\ 5\\ 7 \end{ytableau}\cup \begin{ytableau} 4\\ 5\\ 7 \end{ytableau}\cup \begin{ytableau} 6\\ 8\\ 9 \end{ytableau} },
\end{align*}

\begin{align*}
& X_{ T_{\m}^{(3)}, T_{\n}^{(3)} } = \{(1, 3), (2, 5), (2, 6) \}, \quad 
 Y_{ T_{\m}^{(3)}, T_{\n}^{(3)} } =  \{ (1, 5), (1, 6), (2, 7) \}.
\end{align*}
There is a matching function from $X_{ T_{\m}^{(3)}, T_{\n}^{(3)} }$ to $Y_{ T_{\m}^{(3)}, T_{\n}^{(3)} }$: $(1,3) \mapsto (1,5)$, $(2,5) \mapsto (2,7)$, $(2,6) \mapsto (1,6)$. Therefore ${\rm LC}(T_{\m}^{(3)}, T_{\n}^{(3)})$ is satisfied. 

We have 
\begin{align*}
& X_{ T_{\n}^{(3)}, T_{\m}^{(3)} } = \{(7,1) \}, \quad 
 Y_{ T_{\n}^{(3)}, T_{\m}^{(3)} } =  \{ (7,2) \}.
\end{align*}
There is a matching function from $X_{ T_{\n}^{(3)}, T_{\m}^{(3)} }$ to $Y_{ T_{\n}^{(3)}, T_{\m}^{(3)} }$: $(7,1) \mapsto (7,2)$. Therefore ${\rm LC}(T_{\n}^{(3)}, T_{\m}^{(3)})$ is satisfied. 

For $k=4$, we have  
\begin{align*}
\ch(T_{\m}^{(4)}) \ch(T_{\n}^{(4)}) & = \ch( \scalemath{0.66}{ \begin{ytableau} 5\\ 6\\ 9\\ 10 \end{ytableau} } ) \ch( \scalemath{0.66}{ \begin{ytableau} 1 & 3 & 6\\ 2 & 4 & 7\\ 5 & 6 & 9\\ 8 & 8 & 10 \end{ytableau} } )\\
&  = \ch( \scalemath{0.66}{ \begin{ytableau} 1 & 3 & 5 & 6\\ 2 & 4 & 6 & 7\\ 5 & 6 & 9 & 9\\ 8 & 8 & 10 & 10 \end{ytableau} } ) + \ch( \scalemath{0.66}{ \begin{ytableau} 1 & 6 & 6\\ 2 & 7 & 8\\ 5 &  9 & 9\\  8 & 10 & 10 \end{ytableau} } ) \ch( \scalemath{0.66}{ \begin{ytableau} 3\\ 4\\ 5\\ 6 \end{ytableau} } ) \ne \ch( T_{\m}^{(4)} \cup T_{\n}^{(4)} ).
\end{align*}  
Note that $\scalemath{0.66}{ \begin{ytableau} 3\\ 4\\ 5\\ 6 \end{ytableau} }$ corresponds to the trivial module.

Now we check the conditions ${\rm LC}(T_{\m}^{(4)}, T_{\n}^{(4)})$ and ${\rm LC}(T_{\n}^{(4)}, T_{\m}^{(4)})$. We have 
\begin{align*}
& \scalemath{0.66}{ T_{\m}^{(4)} \sim  \begin{ytableau} 5\\ 6\\ 8\\ 9 \end{ytableau}\cup \begin{ytableau} 6\\ 7\\ 9\\ 10 \end{ytableau}, \quad T_{\n}^{(4)} \sim  \begin{ytableau} 1\\ 2\\ 4\\ 5 \end{ytableau}\cup \begin{ytableau} 2\\ 3\\ 5\\ 6 \end{ytableau}\cup \begin{ytableau} 3\\ 4\\ 5\\ 7 \end{ytableau}\cup \begin{ytableau} 3\\ 4\\ 6\\ 7 \end{ytableau}\cup \begin{ytableau} 4\\ 5\\ 6\\ 8 \end{ytableau}\cup \begin{ytableau} 4\\ 5\\ 6\\ 8 \end{ytableau}\cup \begin{ytableau} 6\\ 7\\ 9\\ 10 \end{ytableau} },
\end{align*}

\begin{align*}
& X_{ T_{\m}^{(4)}, T_{\n}^{(4)} } = \{(1, 3), (1,4), (2, 5), (2, 6) \}, \quad 
 Y_{ T_{\m}^{(4)}, T_{\n}^{(4)} } =  \{ (1, 5), (1, 6), (2, 7) \}.
\end{align*}
There is no matching function from $X_{ T_{\m}^{(4)}, T_{\n}^{(4)} }$ to $Y_{ T_{\m}^{(4)}, T_{\n}^{(4)} }$. Therefore ${\rm LC}(T_{\m}^{(4)}, T_{\n}^{(4)})$ is not satisfied. 
 
We have 
\begin{align*}
& X_{ T_{\n}^{(4)}, T_{\m}^{(4)} } = \{(7,1) \}, \quad 
 Y_{ T_{\n}^{(4)}, T_{\m}^{(4)} } =  \{ (7,1), (7,2) \}.
\end{align*}
There is a matching function from $X_{ T_{\n}^{(4)}, T_{\m}^{(4)} }$ to $Y_{ T_{\n}^{(4)}, T_{\m}^{(4)} }$: $(7,1) \mapsto (7,2)$. Therefore ${\rm LC}(T_{\n}^{(4)}, T_{\m}^{(4)})$ is satisfied.

For $k \ge 5$, we have  
\begin{align*}
& X_{ T_{\m}^{(k)}, T_{\n}^{(k)} } = \{(1, 3), (1,4), (2, 5), (2, 6) \}, \quad 
 Y_{ T_{\m}^{(k)}, T_{\n}^{(k)} } =  \{ (1, 5), (1, 6), (2, 7) \}.
\end{align*}
There is no matching function from $X_{ T_{\m}^{(k)}, T_{\n}^{(k)} }$ to $Y_{ T_{\m}^{(k)}, T_{\n}^{(k)} }$. Therefore ${\rm LC}(T_{\m}^{(k)}, T_{\n}^{(k)})$ is not satisfied. 

We have 
\begin{align*}
& X_{ T_{\n}^{(k)}, T_{\m}^{(k)} } = \{(7,1) \}, \quad 
 Y_{ T_{\n}^{(k)}, T_{\m}^{(k)} } =  \{ (7,1), (7,2) \}.
\end{align*}
There is a matching function from $X_{ T_{\n}^{(k)}, T_{\m}^{(k)} }$ to $Y_{ T_{\n}^{(k)}, T_{\m}^{(k)} }$: $(7,1) \mapsto (7,2)$. Therefore ${\rm LC}(T_{\n}^{(k)}, T_{\m}^{(k)})$ is satisfied.

\end{example}

\section*{}
The authors declare that they have no conflict of interest.

\bibliographystyle{alpha}

\bibliography{biblio}

\end{document}